\documentclass[10pt]{scrartcl} 

\usepackage[utf8]{inputenc}
\usepackage[T1]{fontenc}
\usepackage{lmodern}

\usepackage{amsfonts,amssymb,amsmath,amsthm}
\usepackage{graphicx}  
\usepackage[dvipsnames]{xcolor}
\usepackage{natbib}

\usepackage{authblk}
\usepackage[]{enumitem}
\usepackage{tabularx}
\usepackage{subfig}
\usepackage{mathtools}
\usepackage{url}
\usepackage[ruled,lined]{algorithm2e}
\usepackage{algpseudocode} 
\usepackage{pgfplotstable}
\usepackage{csvsimple}
\usepackage{array,multirow,makecell}

\usepackage[color=blue!30, textsize=scriptsize]{todonotes}
\usepackage{changes}
\definechangesauthor[color=blue!50!black]{MS}
\definechangesauthor[color=green!50!black]{JB}

\DeclareMathOperator*{\cl}{cl}
\DeclareMathOperator*{\conv}{conv}

\newcommand{\WS}{\mathit{WS}} 
\newcommand{\LBS}{\mathcal{L}}
\newcommand{\UBS}{\mathcal{U}}

\newcommand{\XE}{\ensuremath{X_{E}}}
\newcommand{\YN}{\ensuremath{Y_{N}}}

\newcommand{\R}{\mathbb{R}}

\newcolumntype{C}[1]{>{\centering\arraybackslash }b{#1}}
\newcolumntype{L}[1]{>{\arraybackslash }b{#1}}
\newcolumntype{R}[1]{>{\raggedleft \arraybackslash }b{#1}}
\usepackage{tikz, pgfplots}
\usetikzlibrary{arrows, shapes.misc, chains, patterns}
\usetikzlibrary{decorations.pathreplacing}
\usepackage{xifthen}

\tikzset{
  candidat/.style={rectangle, inner sep=0pt, minimum size=0.1cm, draw=gray, fill=gray},
  nds/.style={circle, inner sep=0pt, minimum size=0.12cm, draw=black, fill=black},
  ndns/.style={thick, draw=red, cross out, inner sep=0pt, minimum width=4pt, minimum height=4pt},
  ndns/.style={rectangle, inner sep=0pt, minimum size=0.12cm, draw=black, fill=black},
  test/.style={circle, inner sep=0pt, minimum size=0.12cm, draw=black, fill=black},
}
\pgfplotsset{compat=1.8}

\definecolor{PineGreen}{rgb}{0.0, 0.47, 0.44}
\definecolor{Maroon}{rgb}{0.5, 0.0, 0.0}
\definecolor{RoyalBlue}{rgb}{0.25, 0.41, 0.88}
\definecolor{ChromeYellow}{rgb}{1.0, 0.65, 0.0}

\begin{document}

\title{Adaptive Improvements of Multi-Objective Branch and Bound}

\author{Julius Bau{\ss}}
\author{Sophie N.~Parragh}
\author{Michael Stiglmayr\thanks{Corresponding author}}
\affil{University of Wuppertal, School of Mathematics and Natural Sciences, IMACM, Gaußstr.~20, 42119 Wuppertal, Germany\\
\{bauss,stiglmayr\}@uni-wuppertal.de
}
\affil{Johannes Kepler University Linz, Institute of Production and Logistics Management, Altenberger Str.~69, 4040 Linz, Austria\\
sophie.parragh@jku.at}

\renewcommand{\Affilfont}{\small}
\date{}
\maketitle

\begin{abstract}
Branch and bound methods which are based on the principle ``divide and conquer'' are a well established solution approach in single-objective integer programming. In multi-objective optimization branch and bound algorithms are increasingly attracting interest. However, the larger number of objectives raises additional difficulties for implicit enumeration approaches like branch and bound. Since bounding and pruning is considerably weaker in multiple objectives, many branches have to be (partially) searched and may not be pruned directly. The adaptive use of objective space information can guide the search in promising directions to determine a good approximation of the Pareto front already in early stages of the algorithm. In particular we focus in this article on improving the branching and queuing of subproblems and the handling of lower bound sets.

In our numerical test we evaluate the impact of the proposed methods in comparison to a standard implementation of multiobjective branch and bound on knapsack problems, generalized assignment problems and (un)capacitated facility location problems.


\end{abstract}

\section{Introduction}\label{sec:def}
Multi-objective optimization methods can be generally categorized into objective and decision space algorithms. While objective space methods rely on solving a sequence of scalarized subproblems, decision space algorithms build up solutions in decision space and filter them for dominance like, e.g., multi-objective dynamic programming or multi-objective branch and bound. 
One of the major advantages of decision space methods is that the structure of the feasible set is not investigated several times from scratch.
Objective space methods, however, can make use of commercial single criterion solvers, and despite of the increasing research on  multi-objective branch and bound algorithms within the last decade, some difficulties and problems remain. Since bounding and pruning get considerably weaker with an increasing number of objective functions, many branches have to be searched in-depth and the branch and bound tree can become very large. To avoid the investigation of dominated branches it is important not only to find efficient solutions as soon as possible but also to find a set of (efficient) solutions whose images are well distributed along the
non-dominated frontier. The adaptive use of objective space information can be used to guide the search in promising directions to determine a good approximation of the Pareto front already in early stages of the algorithm.

In this article we consider the solution of multi-objective binary linear optimization problems using multi-objective branch and bound. Our main contribution is thereby the incorporation of objective space information into the decision space search of branch and bound algorithms. In particular, we use an indicator-based scheduling of subproblems and solve adaptively scalarizations to integer optimality in order to improve upper and lower bound sets.

\bigskip
A \emph{multi-objective integer linear program} can be written in the following form:
\begin{equation}\label{eq:MOILP} \tag{MOILP}
   \begin{array}{rr@{\extracolsep{1ex}}c@{\extracolsep{1ex}}ll}
      \min  & \multicolumn{3}{@{\extracolsep{0.75ex}}l}{\displaystyle \bigl( z_1(x),\ldots, z_p(x)\bigr)^\top}\\
      \mathrm{ s.t.} & \displaystyle A\,x &\leq&  b  \\
      &x &\geq& 0 \\
      &x &\in& \mathbb{Z}^n.
   \end{array}
\end{equation}
Since we consider $p \geq 2$ linear objective functions, the objective vector can also be denoted by \(z(x)\coloneqq(z_1(x),\ldots, z_p(x))^\top = C\cdot x \in\R^p\), where $C\in \R^{p\times n}$ is the matrix of objective coefficients. A solution $x\in \R^n$ is called \emph{feasible}, iff $x\in X\coloneqq \{x \in\mathbb{Z}^n: A\,x \leq b, x\geq 0\}$. Hence, $X$ is subset of the so-called \emph{decision space} $\R^n$. Further, the corresponding image $Y\coloneqq \{C\, x\colon x\in X\}$ is a subset of the so-called \emph{objective space} $\R^p$.

In this work we restrict ourselves to binary variables $x\in \{0,1 \}^n$ and therefore to \emph{multi-objective binary linear programs}:

\begin{equation}\label{eq:MO01LP}\tag{MO01LP}
\begin{array}{rr@{\extracolsep{0.75ex}}c@{\extracolsep{0.75ex}}l}
 \min  & \displaystyle z(x) &=& \bigl(z_1(x),\dots,z_p(x)\bigr)^\top\\
 \text{ s.t.} & \displaystyle A\,x &\leq& b  \\
              &       x &\in& \{ 0,1 \}^n.
\end{array}
\end{equation}
Nevertheless, all ideas and approaches that we present in the remaining paper can be easily transferred to \ref{eq:MOILP}.

We use the \emph{Pareto concept of optimality}, since the objective function $z(x)$ is vector-valuated. This concept is based on the componentwise order. Let $y^1,y^2\in \R^p$, then the corresponding dominance relations are given by:
\begin{itemize}
\item $y^1$ \emph{weakly dominates} $y^2$ $(y^1 \leqq y^2)$, if $y^1_k \leq y^2_k$ for $k=1,\ldots ,p$,
\item $y^1$ \emph{strictly dominates} $y^2$ $(y^1 < y^2)$, if $y^1_k < y^2_k$ for $k=1,\ldots ,p$,
\item $y^1$ \emph{dominates} $y^2$ $(y^1 \leq y^2)$, if $y^1 \leq  y^2$ and $y^1 \neq y^2$.
\end{itemize}
We call a feasible solution $x\in X$ \emph{efficient} if there is no other solution $\hat{x} \in X$ with $z(\hat{x}) \leq z(x)$, i.e., there is no other feasible solution that dominates it. The set of all efficient solutions is denoted by $\XE$.  Its corresponding image in the objective space is denoted by $Y_N \coloneqq \{z(x) \in \R^p: x \in \XE \}$. A point $z(x) \in Y_N$ is called \emph{non-dominated}. The description of the set of \emph{weakly efficient} solutions $X_{WE}$ and its image, i.e. the set of \emph{weakly non-dominated} points $Y_{WN}$, can be done analogously. Furthermore, we denote the set of non-dominated points of any set $\mathbb{Q} \subseteq \R^p$ by $\mathbb{Q}_N$. 

The set of non-dominated points $Y_N$ can be decomposed into two subsets, namely the set of \emph{supported non-dominated} points and the set of the \emph{unsupported non-dominated} points. The set of supported non-dominated points is defined by $Y_{SN} \coloneqq \{z(x)\in Y_N: z(x) \in (\conv(Y)+\R^p_{\geqq})_N \}$, where $\conv(Y)$ is the convex hull of $Y$ and $\R^p_{\geqq} \coloneqq \{y \in \R^p:y\geqq 0\}$. Thus, the set of unsupported non-dominated points is defined by $Y_{UN} \coloneqq \{z(x) \in Y_N : z(x) \notin Y_{SN} \}$. The supported points are located on the boundary of the convex hull of $Y$ and the unsupported points are located in its (relative) interior. 

As the solution of a \ref{eq:MOILP} we consider the computation of all non-dominated points $\YN$ and a corresponding \emph{minimal complete set} of efficient solutions. More precisely, it is required to find all non-dominated points $y\in Y_N$ and for each non-dominated point a corresponding solution $\hat{x}$ with $z(\hat{x}) = y \in Y_N$.

As its name implies, one of the key components of a branch and bound algorithm is the bounding of the set of non-dominated points. Although bounding is straightforward and convenient in the single objective case, the computation of lower and upper bounds is crucial in the multi-objective case. The \emph{ideal point} $y^I$ and the \emph{Nadir point} $y^N$ are well known and also the tightest componentwise lower and upper bounds of the set of non-dominated points $Y_N$. They are defined as follows:

\[
  y^I_k = \min_{y\in Y} y_k  \quad \text{and}\quad y^N_k = \max_{y\in Y_N} y_k \qquad \text{for } k=1,\ldots p. 
\]
Obviously, those points are valid bounds of $Y_N$, since $y^I \leqq y \leqq y^N$ for all $y \in Y_N$.  Although these bounds are a direct extension of the single-objective case, they are rather weak. Since there is in general no optimal solution which optimizes all objectives at once, it holds $y^I \neq y^N$ in most of the cases. As a result, instead of using a single point as bound, \emph{bound sets} are used. We refer to \citep{ehrgott2007bound} for the following definitions of bound sets.
\begin{itemize}
   \item A \emph{lower bound set}  $\LBS \subset \R^p$ for $Y_N$ is a $\R^p_\geqq$-closed (i.e., the set \(\LBS +\R^p_\geqq\) is closed), $\R^p_\geqq$-bounded (i.e., there exists a \(y\in\R^p\) such that \(\LBS\subset y+\R^p_\geqq\)) stable set (i.e., $\LBS \subset (\LBS +\R^p_\geqq)_N$), such that $Y_N \subset (\LBS +\R^p_\geqq)$.
   \item An \emph{upper bound set} $\UBS \subset \R^p$ for $Y_N$ is a $\R^p_\geqq$-closed, $\R^p_\geqq$-bounded, stable set,      such that $Y_N \subset \cl\bigl((\UBS +\R^p_\geqq)^\complement \bigr)$.
\end{itemize}
The lower bound sets and upper bound set that we will use in the following sections for our branch and bound methods will fulfill all those requirements.

The remainder of this paper is structured as follows: In Section~\ref{sec:relwork} we summarize the components of multi-objective branch and bound algorithms and review the literature on multi-objective branch and bound. Our augmentations of multi-objective branch and bound using objective space information are presented and discussed in Section~\ref{sec:aug}.
The numerical results in Section~\ref{sec:res} show the effectiveness of the proposed improvements on knapsack, generalized assignments and facility location problems. 
Section~\ref{sec:dis} concludes this article.

\section{Related Work: Multi-objective Branch and Bound Frameworks in the Literature}\label{sec:relwork}
The majority of objective space methods, for example the weighted sum method or the $\varepsilon$-constraint method, solve a scalarized integer program from scratch in each iteration without transfering starting solutions found in prior iterations. Hence, an increasing research interest in decision space methods, mainly the branch and bound method, can be observed in the last decades. 

A multi-objective branch and bound operates in the same way as its well-known single-objective version. Since the problem, which needs to be solved, is too hard to be solved directly, it is divided into easier subproblems. Every created subproblem is associated with a node, resulting in a tree data structure. Thereby, a node $i$ is the child node of node $j$ if and only if the fesaible set of the subproblem corresponding to node $i$ is a subset of the feasible set of the (sub)problem corresponding to node $j$. One of the first, if not the first, multi-objective branch and bound with an underlying tree structure was developed by \cite{Klein1982an}. In each iteration an active node is selected and its corresponding lower bound set is computed. We start, obviously, with the root to which the original problem is associated. After the computation of the lower bound set, we possibly update the upper bound set and check if it is possible to fathom the node. If the node cannot be fathomed, the corresponding problem has to be further divided into new subproblems (branching). As a result a branch and bound method is made up of the following components: node selection, lower bound, upper bound, fathoming and branching. Since the choice of each component is crucial for the performance of this method, we are going to present some of the most frequently used approaches of each component.

\paragraph{Lower bound:}
 In a branch and bound approach, in each iteration a node is selected and its lower bound is computed. Since bounding is considerably weaker in the multi-objective case, the choice of the bound computation approach is crucial for the performance of those methods.
 Hence, different approaches for the computation of the lower bound have been proposed in the last decades.

One of the first approaches was the so-called \emph{minimal completion} \citep[see, e.g.,][]{Klein1982an,Kiziltan1983an}. In this approach each variable was fixed to $0$ or $1$ (for \ref{eq:MO01LP}) depending on the corresponding objective values. Although such a solution is integer it is not necessarily feasible, since it could violate certain constraints. Another approach is to use the ideal point as the lower bound. Due to its expensive computation it is often replaced by using the ideal point of the \emph{linear relaxation} \citep[see, e.g.,][]{Mavrotas1998a,Mavrotas2005multi}.

\citet{Sourd2008a} proposed to use the solution of the \emph{convex relaxation} as lower bound set \citep[see also][]{Vincent2013multiple}. This convex relaxations can be solved by, e.g., using a dichotomic scheme \citep[]{aneja1979bicriteria,oezpeynirci_2010an,Przybylski2010a}.

Like in the single-objective case, a lower bound set is most commonly obtained by solving the linear relaxation. It is used, for example, by \citet{Vincent2013multiple}, \cite{Adelgren2022branch} and \cite{Parragh2019branch}. Besides using the dichotomic schemes it can also be obtained with \emph{Benson's outer approximation algorithm} \citep{Benson1998an,ehrgott12dual}. The algorithm is initialized with a (weak) lower bound which is improved in each iteration. The outer approximation approach ensures that the algorithm can be aborted at any time since the produced lower bound set is valid all the time. In \citet{Forget2022warm} an approach for warm-starting this outer approximation algorithm is proposed.

Note that we denote a lower bound set $\LBS$ as \emph{convex lower bound set} if the set $\LBS+ \R^p_\geq$ is convex.

\paragraph{Upper bound:}
In the single-objective case the usual upper bound is given by the best found solution so far. In the multi-objective case, the direct execution of this idea is used. A bound set $\UBS$ is stored in a so-called \emph{incumbent list}. During the algorithm all feasible solutions and their corresponding images in the objective space are stored in this list if they are not dominated by another feasible solution found so far. After computing the lower bound set, its extreme points are checked for integer feasibility. Such a feasible solution $\hat{x} \in X$ is then added to $\UBS$ if there is no other $\bar{x} \in X$ such that $z(\bar{x}) \leq z(\hat{x})$. All solutions $x\in \UBS$ that are dominated by a newly added solution $\hat{x}$ are removed from the upper bound set. 
\[
\UBS \coloneqq 
\begin{cases}
	\UBS & \text{if}\; \exists x\in\UBS\colon C(x)\leq C(\hat{x}) \\
	\{\hat{x}\}\cup \{x\in \UBS \colon C(\hat{x}) \nleq C(x) \} & \text{otherwise}.
\end{cases}
\]

\paragraph{Node selection:} 
In each iteration of a branch and bound method, the first task is to select an unexplored node from the tree of subproblems. We call this selected node \emph{active}. Since the order of the considered nodes have a significant impact on the number of created nodes and the computational time, several approaches have been proposed in the literature. We distinguish between \emph{static strategies} and \emph{dynamic strategies}. While dynamic strategies consider information gained in previous iterations for the choice of the next node, static strategies consider the nodes in a consistent order, e.g., \emph{first in first out} or \emph{last in first out}.

Most of the branch and bound methods proposed in the literature use static strategies. The best-known static strategies are the \emph{depth-first strategy} and the \emph{breadth-first strategy}. The depth-first strategy is used, for example, by \citet{Ulungu1997solving}, \citet{Kiziltan1983an}, \citet{Visee1998two} and \citet{Vincent2013multiple}. In, for example, \citet{Parragh2019branch} the breadth-first strategy is used after showing that it is more convenient for some of their problem classes. 

Although it is common to use dynamic strategies in single-objective branch and bound methods they are rarely applied in the multi-objective case. Nevertheless several different dynamic approaches have been proposed in the last decades. One of them is to rely on the optimal objective value of a linear relaxation of a weighted sum scalarization. In the first place this  relaxation is solved to generate the lower bound set of the corresponding subproblem. This technique is used by, e.g., \citet{Stidsen2014a} and \cite{Gadegaard2019bi}. In \citet{Belotti2013a} another dynamic strategy is proposed. For all solutions computed by the linear relaxation (lower bound computation) the integer infeasibility  is summed up for all (integer) variables. Then, the node with the largest sum of integer infeasibility is selected. One of the most used selection rules in the single-objective case is to choose the node with the largest gap between the lower and upper bound. But since the bounds are given by sets with multiple points in the multi-objective case, there are numerous options to measure the gap. In \citet{Jesus2021on} the gap is computed with the so-called $\varepsilon$-indicator. \citet{Adelgren2022branch} use a slightly adapted Hausdorff distance as gap measure. \citet{Forget2023enhancing} compare this gap measurement and a best bound value derived from solving a weighted sum problem. In \citet{Bauss2023augmenting} the gap is obtained by computing the so-called \emph{approximated hypervolume gap}.

\paragraph{Fathoming:}
In a worst-case scenario a branch and bound approach produces the total enumeration of all feasible solutions. To avoid this, there are rules that allow us to fathom a node, since there are no new non-dominated points in the corresponding branch. The three different cases that might occur are: infeasibility, optimality or dominance. 
\begin{itemize}
	\item[i) ] \emph{Fathoming by infeasibility}: This case is a direct extension of the single-objective branch and bound. When solving a relaxed problem to compute the lower bound set returns infeasibility, then the corresponding subproblem is infeasible as well. This is obviously true since the feasible set of the subproblem is a subset of the feasible set of its relaxation.
	\item[ii) ] \emph{Fathoming by optimality:} Like in the single-objective case a node can be fathomed by optimality if the upper bound $\UBS$ is equal to the lower bound $\LBS$. This would imply that this node is solved to optimality and that there is no need to divide it further. But since the bound sets in general consists of multiple points this rarely happens. The only possibility where this rule can be applied is when the lower and upper bound set consist of the same single point, i.e., the ideal point.
	\item[iii) ] \emph{Fathoming by dominance:} A node can be fathomed by dominance if all points of the lower bound set are dominated by at least one point of the upper bound set, i.e., if all feasible solutions of this subproblem are dominated by points of the incumbent list. This dominance check might lead to difficulties since the bound sets could possibly have different properties. Although the upper bound set, i.e., the incumbent list, is always a finite set of points for (\ref{eq:MO01LP}), the lower bound can be computed in different ways. Therefore, the lower bound set can be composed of an unique point (e.g., the ideal point), a unique hyperplane (e.g., obtained by solving a single weighted sum scalarization) or a complex lower bound (e.g., obtained by solving the linear relaxation). Depending on the shape of the lower bound, there are different approaches to verify dominance. Since a large number of the more recently published papers rely on solving the convex relaxation, solving the linear relaxation or using a unique hyperplane for the lower bound set, we refer to the dominance test proposed in \citet{Sourd2008a}. For those lower bound sets it is sufficient to check if all local upper bounds are located below the lower bound set \citep[see][]{Sourd2008a,Gadegaard2019bi}. An overview of existing fathoming by dominance rules can be found in \citet{Belotti2016fathoming}.
\end{itemize}

\paragraph{Branching:} 
If a node cannot be fathomed its corresponding problem is further divided into smaller subproblems. Since we are discussing (\ref{eq:MO01LP}) two subproblems are created by fixing a variable to $0$ or respectively to $1$. This procedure results in a binary tree structure. 

Of course the choice of the variable on which the branching is applied is crucial for the performance of the branch and bound algorithm. Therefore several branching rules have been proposed. Again, we need to distinguish between static and dynamic strategies. If the order of the variables that are branched is known in advance, we call this strategy static. In every iteration, the first variable of the list which is not fixed in the corresponding subproblem is selected. The most basic idea for single-objective knapsack problems  is to choose the most promising variable according to the \emph{profit-to-weigth ratio} $c_i \slash a_i, i \in \{1,\dots,n\}$, where $c\in \R^n$ is the objective vector and $a\in \R^n$ is the vector of constraint coefficients \citep{Kellerer2004knapsack}. Although there is no direct extension to the multi-objective case, there are some approaches that consider the conflicting objective function values \citep[see, e.g.,][]{Bazgan2009solving,Ulungu1997solving}.

Dynamic branching strategies take information into account that have been gained throughout the algorithm and base their variable selection on it. Although most of the published papers use static strategies there are several dynamic approaches. One of them is to count how often a variable is not integer in all solutions corresponding to the extreme points of the computed lower bound set and then to choose the variable, which is most often fractional. Another approach is to sum up the integer infeasibility and choose the variable that is most fractional \citep{Belotti2013a}. \citet{Stidsen2014a}, \cite{Stidsen2018a} and \citet{Gadegaard2019bi} give a single-objective solver the choice of their next variable to branch on.

\paragraph{}Additionally to those crucial components of a branch and bound method there are further approaches that contain more features. One of the more recent approaches is to also include \emph{objective space branching}. This procedure subdivides problems by adding constraints that touch the objective functions. This idea and the similar idea of \emph{Pareto branching} are proposed in \citet{Stidsen2014a}. Pareto branching creates additional subproblems by adding upper bounds on the objective functions to disregard areas in the objective space, where the lower bound is already dominated. It is also used by, e.g., \citet{Parragh2019branch}, \citet{Gadegaard2019bi} and \citet{Adelgren2022branch}.

In \citet{Bauss2023augmenting} improvements for the lower and upper bound set have been proposed. By adaptively solving scalarizations to integer optimality, objective space information are gained, which can be used to update the upper bound set and improve the lower bound set in possibly all created subproblems. This bound improvements lead to a higher fathoming rate.

A survey on multi-objective branch and bound approaches is provided in \citet{przybylski17multi}, a survey for its single-objective counterpart can be found in \citet{Morrison2016branch}.

\section{Using Objective Space Information in Multi-objective Branch and Bound}\label{sec:aug}
In this section we show how objective space information can improve the computational performance of multi-objective branch and bound. Thereby we focus on dynamic branching strategies and the usage of objective space information to improve the computational efficiency of multi-objective branch and bound approaches. 

As already mentioned, there are two main shortcomings of multi-objective branch and bound algorithms. Firstly, the bounding is weaker, compared to its single-objective counterpart and secondly, optimized single-objective solvers lead to the supremacy of objective space methods. Therefore, we utilize those single-objective solvers to solve scalarized integer programs and use the obtained information to possibly improve the lower and upper bound set.
Furthermore, we present dynamic node selection strategies. Although it is common to use dynamic strategies in the single-objective case, they are rarely applied in the multi-objective case. 

In order to present the improved components of our multi-objectve branch and bound its basic structure is summarized in Algorithm~\ref{alg:BB}.
\begin{algorithm}
	\caption{Multi-objective Branch and Bound Algorithm}\label{alg:BB}
	\begin{description}[font=\rmfamily\itshape]
    \item[Step 0] Let $\nu^0$ be the node corresponding to the initial problem. 
      Initialize the list of nodes $\mathcal{N}\coloneqq \{\nu^0\}$ and the upper bound set $\UBS\coloneqq \emptyset$.
    \item[Step 1] Select and remove a node $\nu$ from $\mathcal{N}$.
    \item[Step 2] Compute the lower bound set for the (sub)problem corresponding to $\nu$.
    \item[Step 3] If a new integer feasible solution has been found update $\UBS$ if necessary. 
    \item[Step 4] Check whether node $\nu$ can be fathomed. If $\nu$ can be fathomed go to Step 1.
    \item[Step 5] Create two disjoint subproblems. Add the corresponding newly created nodes to $\mathcal{N}$. If $\mathcal{N} \neq \emptyset$ go to Step 1.
	\end{description}
\end{algorithm}

\subsection{A Dynamic Node Selection Strategy} \label{subsec:nodesel}
The order in which nodes in a branch and bound tree are explored has a significant impact on the total number of nodes and therefore on the computation time. Choosing the right nodes and obtaining good feasible solutions in early stages of the algorithm can lead to a decrease of nodes to explore, since they could dominate numerous other nodes, which as a result can be fathomed. The arising difficulty is to use a good node selection strategy that causes a desirable reduction of nodes. Therefore, instead of using a static strategy, like, e.g., the depth-first strategy, we propose a new dynamic strategy.

The underlying principle of this strategy is a direct extension of a node selection strategy that is frequently used in the single-objective case. There, in each iteration the node with the largest gap between upper and lower bound is selected \citep[see, e.g.,][]{Dechter1985generalized, Morrison2016branch}. In the multi-objective case there are numerous approaches to measure the gap between the lower bound set and the set of local upper bounds and/or the points in the incumbent list. See \cite{bauss24adapting} for a comparison of different gap measures and the according queuing of subproblems. We propose to use the so-called \emph{approximate hypervolume gap} to select the active subproblem from the queue, where we rely on the definition of hypervolume of \citet{Zitzler1999multiobjective}. Like in its single-objective counter part, in every iteration the node with the largest gap is selected. Thereby, we distinguish two different ways to approximate the hypervolume gap.

The first approach is called the \emph{local hypervolume gap} \citep[see][for the bi-objective case]{Bauss2023augmenting}. Instead of measuring the total gap between the lower and upper bound set, we only consider the largest gap between a single local upper bound and the lower bound set, i.e., we only consider the volume of the largest search zone. In \citet{Klamroth2015on} a detailed analysis of search zones, search regions, corresponding local upper bounds and their defining nondominated points is given.

\begin{figure}[htbp!]
	\subfloat[\label{hyperv-a}]{
			\begin{tikzpicture}[scale=0.5]
			\draw[style=help lines] (0,0) grid (12,12);
			\draw[line width=1.5,->,black!75](0,0) -- (12,0) node[right] {$z_1$};
			\draw[line width=1.5,->,black!75] (0,0) -- (0,12) node[above] {$z_2$};
			
				\fill[lightgray] (1.15,9) -- (6,9) -- (6,1.1);

				\draw (2,9) node[above] {$z^1$};
			\draw (6,7) node[left] {$z^2$};
\draw (9,5) node[left] {$z^3$};
\draw (10,1) node[left] {$z^4$};
\draw (6,9) node[above right] {$lu^1$};
\draw (9,7) node[above right] {$lu^2$};
\draw (10,5) node[above right] {$lu^3$};
			\draw[RoyalBlue,line width = 1.2] (1,12)--(1,10.5) -- (1.5,5) -- (3,2) -- (8,0.5)-- (12,0.5);
			\draw[dashed] (2,9) -- (6,9) -- (6,7) -- (9,7) -- (9,5) -- (10,5) -- (10,1);			
			\filldraw(2,9) circle (3pt);
			\filldraw  (10,1) circle (3pt);
			
			\filldraw (9,5) circle (3pt);
			\filldraw (6,7) circle (3pt);
			
			\filldraw[Maroon] ([xshift=-3pt,yshift=-3pt]6,9) rectangle ++(6pt,6pt);
			\filldraw[Maroon] ([xshift=-3pt,yshift=-3pt]9,7) rectangle ++(6pt,6pt);
			\filldraw[Maroon] ([xshift=-3pt,yshift=-3pt]10,5) rectangle ++(6pt,6pt);
			
	\end{tikzpicture}}	
	\subfloat[\label{hyperv-b}]{
			\begin{tikzpicture}[scale=0.5]
			\draw[style=help lines] (0,0) grid (12,12);
			\draw[line width=1.5,->,black!75](0,0) -- (12,0) node[right] {$z_1$};
			\draw[line width=1.5,->,black!75] (0,0) -- (0,12) node[above] {$z_2$};
			
			\fill[lightgray] (9,7) -- (9,0.5) -- (1.3,7);

			\draw (2,9) node[above] {$z^1$};
			\draw (6,7) node[left] {$z^2$};
			\draw (9,5) node[left] {$z^3$};
			\draw (10,1) node[left] {$z^4$};
			\draw (6,9) node[above right] {$lu^1$};
			\draw (9,7) node[above right] {$lu^2$};
			\draw (10,5) node[above right] {$lu^3$};
			\draw[RoyalBlue,line width = 1.2] (1,12)--(1,10.5) -- (1.5,5) -- (3,2) -- (8,0.5)-- (12,0.5);
			\draw[dashed] (2,9) -- (6,9) -- (6,7) -- (9,7) -- (9,5) -- (10,5) -- (10,1);			
			\filldraw(2,9) circle (3pt);
			\filldraw  (10,1) circle (3pt);
			
			\filldraw (9,5) circle (3pt);
			\filldraw (6,7) circle (3pt);
			
			\filldraw[Maroon] ([xshift=-3pt,yshift=-3pt]6,9) rectangle ++(6pt,6pt);
			\filldraw[Maroon] ([xshift=-3pt,yshift=-3pt]9,7) rectangle ++(6pt,6pt);
			\filldraw[Maroon] ([xshift=-3pt,yshift=-3pt]10,5) rectangle ++(6pt,6pt);
	\end{tikzpicture}}
	
	\subfloat[\label{hyperv-c}]{
			\begin{tikzpicture}[scale=0.5]
			\draw[style=help lines] (0,0) grid (12,12);
			\draw[line width=1.5,->,black!75](0,0) -- (12,0) node[right] {$z_1$};
			\draw[line width=1.5,->,black!75] (0,0) -- (0,12) node[above] {$z_2$};
			
			\fill[lightgray] (1,9) -- (6,9) -- (6,0.5) -- (1,0.5);

			\draw (2,9) node[above] {$z^1$};
			\draw (6,7) node[left] {$z^2$};
			\draw (9,5) node[left] {$z^3$};
			\draw (10,1) node[left] {$z^4$};
			\draw (6,9) node[above right] {$lu^1$};
			\draw (9,7) node[above right] {$lu^2$};
			\draw (10,5) node[above right] {$lu^3$};
			\draw[RoyalBlue,line width = 1.2] (1,12)--(1,10.5) -- (1.5,5) -- (3,2) -- (8,0.5)-- (12,0.5);
			\draw[dashed] (2,9) -- (6,9) -- (6,7) -- (9,7) -- (9,5) -- (10,5) -- (10,1);			
			\filldraw(2,9) circle (3pt);
			\filldraw  (10,1) circle (3pt);
			
			\filldraw (9,5) circle (3pt);
			\filldraw (6,7) circle (3pt);
			
			\filldraw[Maroon] ([xshift=-3pt,yshift=-3pt]6,9) rectangle ++(6pt,6pt);
			\filldraw[Maroon] ([xshift=-3pt,yshift=-3pt]9,7) rectangle ++(6pt,6pt);
			\filldraw[Maroon] ([xshift=-3pt,yshift=-3pt]10,5) rectangle ++(6pt,6pt);
			
					\filldraw[RoyalBlue] ([xshift=-3pt,yshift=-3pt]1,0.5) rectangle ++(6pt,6pt);
			
	\end{tikzpicture}}
\subfloat[\label{hyperv-d}]{
	\begin{tikzpicture}[scale=0.5]
		\draw[style=help lines] (0,0) grid (12,12);
		\draw[line width=1.5,->,black!75](0,0) -- (12,0) node[right] {$z_1$};
		\draw[line width=1.5,->,black!75] (0,0) -- (0,12) node[above] {$z_2$};
		
		\fill[lightgray] (9,7) -- (9,0.5) -- (1,0.5) -- (1,7);

		\draw (2,9) node[above] {$z^1$};
		\draw (6,7) node[left] {$z^2$};
		\draw (9,5) node[left] {$z^3$};
		\draw (10,1) node[left] {$z^4$};
		\draw (6,9) node[above right] {$lu^1$};
		\draw (9,7) node[above right] {$lu^2$};
		\draw (10,5) node[above right] {$lu^3$};
		\draw[RoyalBlue,line width = 1.2] (1,12)--(1,10.5) -- (1.5,5) -- (3,2) -- (8,0.5)-- (12,0.5);
		\draw[dashed] (2,9) -- (6,9) -- (6,7) -- (9,7) -- (9,5) -- (10,5) -- (10,1);			
		\filldraw(2,9) circle (3pt);
		\filldraw  (10,1) circle (3pt);
		
		\filldraw (9,5) circle (3pt);
		\filldraw (6,7) circle (3pt);
		
		\filldraw[Maroon] ([xshift=-3pt,yshift=-3pt]6,9) rectangle ++(6pt,6pt);
		\filldraw[Maroon] ([xshift=-3pt,yshift=-3pt]9,7) rectangle ++(6pt,6pt);
		\filldraw[Maroon] ([xshift=-3pt,yshift=-3pt]10,5) rectangle ++(6pt,6pt);
		
		\filldraw[RoyalBlue] ([xshift=-3pt,yshift=-3pt]1,0.5) rectangle ++(6pt,6pt);
		
\end{tikzpicture}}
	\caption{A bi-objective example of computation of the two different approximated hypervolume gap approaches. In (a) and (b), the the approximated hypervolume gap (gray) is visualized for the local upper bounds $lu^1$ and $lu^2$. In (c) and (d), the hypervolume of the box (gray) defined by the local ideal point and the local upper bound $lu^1$ respectively $lu^2$ is shown.}
	\label{fig:hyperv}
\end{figure}
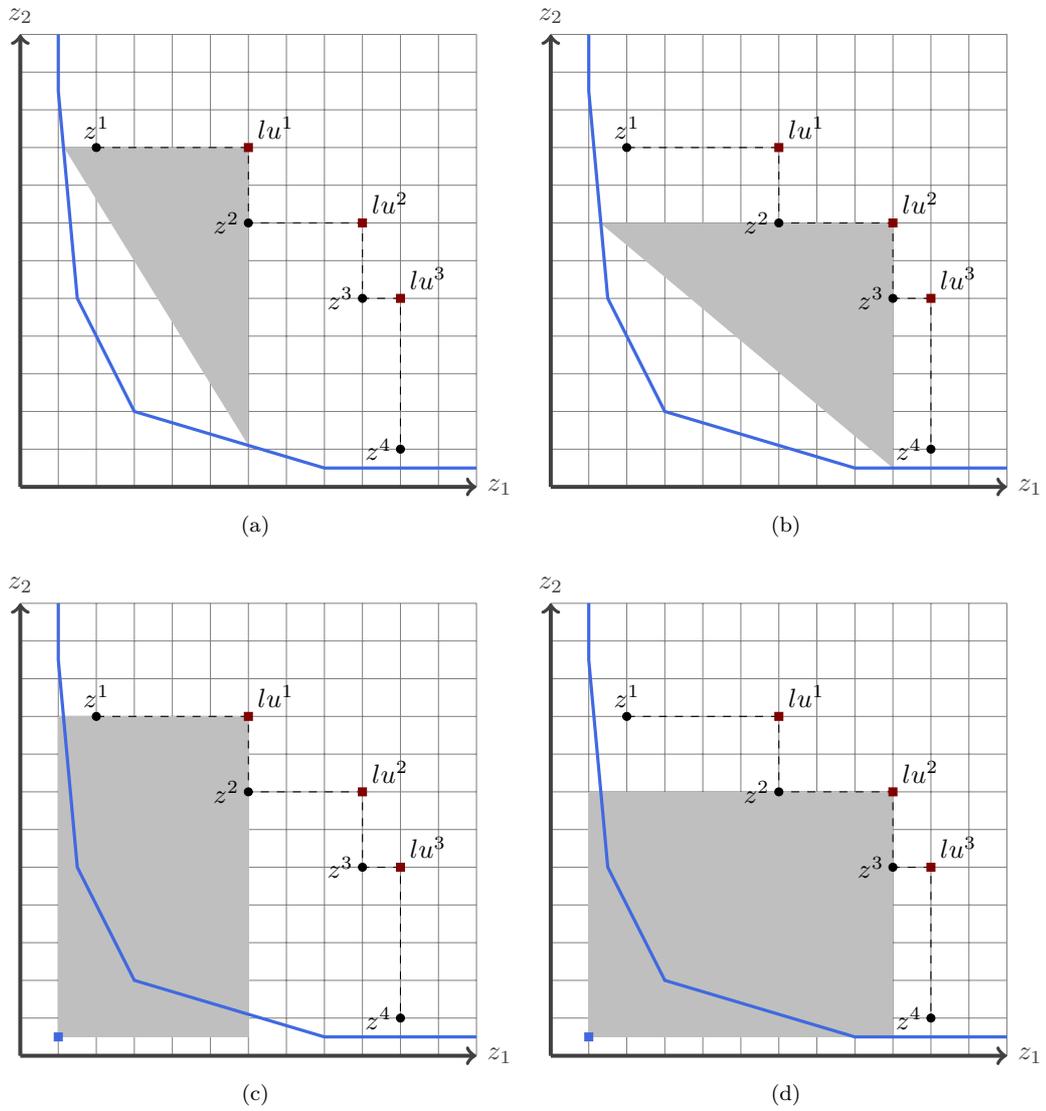

Figures \ref{hyperv-a} and \ref{hyperv-b} illustrate the local hypervolume  gap approach for a bi-objective example.  The points $z^1,\dots ,z^4 \in \UBS$ are points  of the incumbent list $\UBS$. Furthermore, the points $lu^1,\dots,lu^3$ are their corresponding local upper bounds. Note that we only need to consider those local upper bounds that are located above the lower bound set $\LBS$ of the active node $\nu$. The lower bound set $\LBS$ is illustrated by the blue line and is obtained by solving the linear relaxation of the corresponding subproblem in node $\nu$. For every local upper bound the approximated volume of the corresponding search zone is computed as the hypervolume of the simplex spanned by a local upper bound and $p$ spanning points on the lower bound set. Those spanning points w.r.t.\ a local upper bound $lu$ are defined as its axis-parallel projections on the lower bound set, i.e., $sp^i(lu) \coloneqq \{l \in \LBS \colon l_{p+1-i} = lu_{p+1-i}\}$ for  $i =1,\dots,p$. Thus, the hypervolume gap between a local upper bound $lu$ and the lower bound set $\LBS$ is given by
\[
  hg(lu) \coloneqq \frac{|\det(G)|}{p!},
\]
with $G\coloneqq (sp^1-lu,\dots,sp^p-lu)\in \R^{p\times p}$.
Let $K\subset \R^p$ be the set of all local upper bounds based on the incumbent list in node $\nu$. Then, the local hypervolume gap of the node $nu$ is defined as the largest hypervolume spanned by a local upper bound in $K$ and the corresponding spanning points, i.e.,
\[
  lhg(\nu) \coloneqq \max_{i = 1,\dots |K|} hg(lu^i).
\]
Obviously, $lhg(\nu)$ is in general only a rough approximation since the real hypervolume of the search zone is underestimated by neglecting possibly large parts. Even though this approximation simplifies the computation significantly compared to the computation of the real hypervolume of a search zone, it gets too time consuming with an increasing number of objective functions. 
Both the number of local upper bounds and the number of facets of the lower bound set increase substantially with the number of objectives. Since the projection of the local upper bounds on the facets of the lower bound set to determine the spanning points requires a significant amount of the total computation time, we simplify the computation at this point further. 

The second approach to measure the gap between the lower bound $\LBS$ and the upper bound $\UBS$ is to compute the \emph{hypervolume of a search zone box} that is defined by a local upper bound $lu$ and the local ideal point of the lower bound set $l^I$, defined by $l^I_k\coloneqq \min_{l\in \LBS} l_k, k \in \{1,\dots,p\}$. Therefore, the hypervolume of the box defined by local upper bound $lu$ is given by
\[hb(lu)\coloneqq \prod_{k=1}^p (lu_k - l^I_k).\]
Again, this is computed for every local upper bound located above the lower bound and afterwards the volume of the largest box is assigned to the corresponding node. Consequently, the gap between the lower bound and upper bound in node $\nu$, using the hypervolume of a search zone approach, is given by
\[hsz\coloneqq \max_{i=1,\dots,|K|} hb(lu^i).\]

When new nodes are created by branching, the approximated hypervolume gap of the parent node is assigned to the child nodes to avoid the computation of the lower bound set before the child node becomes active. Note that the set of local upper bounds is initialized with the point $(M,\ldots,M)^\top \in \R^p$ with a sufficiently large value $M\gg 0$. This allows us to immediately apply this node selection strategy at the beginning of our algorithm.

Both presented gap measures are illustrated for a bi-objective example in Figure \ref{fig:hyperv}. The local hypervolume gap approach is shown in Figures \ref{hyperv-a} and \ref{hyperv-b}, where the gray area illustrates the hypervolume of the simplex spanned by a local upper upper bound its corresponding spanning points. In Figures \ref{hyperv-c} and \ref{hyperv-d}, the hypervolume of a search zone box approach is illustrated. The gray area indicates the hypervolume of the box spanned by a local upper bound and the local ideal point of the lower bound, which is illustrated by the blue rectangle.

In Algorithm~\ref{alg:BB}, the node with the largest assigned hypervolume gap is selected in Step 1. The value of the hypervolume gap is updated in Step 4 if the node can not be fathomed.

\subsection{Improving Bounds by Solving Scalarizations to Integer Optimality}
In this section, we propose ways to integrate objective space methods into a branch and bound framework. By solving suitable scalarizations to integer optimality, we obtain non-dominated points and thus objective space information that can be used to improve the lower as well as the upper bound set. Let $\hat{x}$ be the integer optimal solution obtained by solving a scalarization of the underlying problem. Since $z(\hat{x})$ is non-dominated, it can be added to the incumbent list (if not contained already), which improves the upper bound set. Additionally, depending on the used scalarization technique, the lower bound set might be improved. An improved lower bound reduces the area where possibly new non-dominated points could be found and an improved upper bound set leads to a higher fathoming rate.

\subsubsection{Warmstarting the Bound Sets}\label{subsub:warmstart}
A branch and bound algorithm benefits from good bound sets and the earlier good bounds are obtained the more impact it has on the performance. Therefore we present a warmstarting approach for the bound set. For this purpose we introduce the weighted sum scalarization that is defined as
\begin{equation}\label{eq:WS}\tag{$\WS_\lambda$}
	\begin{split}
		\min \;\; &  \sum_{i=1}^p \lambda_i\, z_i(x)\\
		\mathrm{s.t.}\;\; & x \in X.
	\end{split}
\end{equation}
Every optimal solution of \ref{eq:WS} is at least weakly efficient for $\lambda \in \R^p_{\geq} \coloneqq \{\lambda\in \R^p\colon \lambda\geq 0\}$  but efficient for $\lambda \in \R^p_{>}$, regarding a (\ref{eq:MO01LP}). Note that this scalarization can only determine supported non-dominated points \citep[see, e.g.,][]{ehrgott05multicriteria}.

We solve a limited number of weighted sum scalarizations with different weights $\lambda \in \R^p_>$ in a preprocessing step of the branch and bound algorithm. This produces a warmstart of the lower bound set and the incumbent list. For the bi-objective case this idea is proposed in \citet{Boekler2021an}. The authors use an outer approximation algorithm to generate $\conv(Y)_N$, that can be used as an initial lower bound set in a multi-objective branch and bound algorithm. However, this approach is way more difficult for $p\geq 3$. It has to deal with similar difficulties as the dichotomic search scheme approach \citep[see][]{aneja1979bicriteria,Przybylski2010a,przybylski2019simple}. Hence, we use a predefined weight set $\Lambda$ to overcome these problems. For every $\lambda \in \Lambda$, a weighted sum scalarization is solved in the root node of the branch and bound tree. 

For each $\lambda \in \Lambda$ the scalarized problem \eqref{eq:WS} can be solved with a single-objective integer linear programming solver. Let $\hat{x}$ be the optimal solution of \eqref{eq:WS}, then $z(\hat{x})$ is a (supported) non-dominated point. Hence, we can add this point to the incumbent list if it is not contained already, and update the list of corresponding local upper bounds. Additionally we obtain further objective space information, that can possibly improve the lower bound set of all nodes that are explored during the algorithm. The optimal solution yields a level set $\{y\in\R^p\colon \lambda^\top y = \lambda^\top z(\hat{x}) \}$, which implies the valid inequality $\lambda^\top z(x) \geq \lambda^\top z(\hat{x})$ for all $x\in X$. Since this inequality is obtained by solving a scalarization of the root node problem, it holds for every subproblem. In \citet{Bauss2023augmenting} the level sets were integrated in an already computed lower bound set by computing the potential cuts of a $\LBS$ and the level set. However, if more than two objectives are considered this would require a lot more computational time, because of the large amount of geometrical operations. Therefore, we avoid this problem by adding the corresponding  inequality to the integer programming formulation of every subproblem. During the algorithm the lower bounds are improving which might cause redundancy w.r.t.\ these inequalities. Hence, in each iteration we check for redundancy in the active node and delete redundant inequalities.

In Algorithm~\ref{alg:BB}, the warmstart of the lower bound set is integrated in Step 0.

\subsubsection{Improving the Upper Bound set by Using $\varepsilon$-constraint Scalarization}
Since the weighted sum approach, which is applied as a warmstarting technique for the lower bound set and the incument list, only generates supported non-dominated points, it might be useful to also use other scalarization techniques, that also compute unsupported non-dominated points. In \citet{Bauss2023augmenting} the \emph{augmented weighted Tchebycheff scalarization} is used to improve the upper bound by obtaining possibly unsupported non-dominated points. Additionally, the corresponding level set is used to improve the lower bound set. Unfortunately, the updated lower bound set had nearly no impact on the performance, since the lower bound improves just locally and is computationally hard to handle due to its (in general) non-convex structure. 

Since we are not updating the lower bound set in this step and only aim at computing unsupported efficient solutions we rely on the \emph{$\varepsilon$-constraint scalarization}, which was firstly introduced by \citet{Haimes1971on}. Thereby, one of the objective functions $z_k,k \in \{1,\dots,p\}$ of  (\ref{eq:MO01LP}) is selected as the objective function of the scalarized problem. The remaining $p-1$ objective functions are transformed into constraints that bound the corresponding objective values. Hence, the $\varepsilon$-constraint scalarization can be written in the form:
\begin{equation}\label{eq:econs} \tag{$\varepsilon$-C}
	\begin{array}{rr@{\extracolsep{1ex}}c@{\extracolsep{1ex}}ll}
		\min  & \multicolumn{3}{@{\extracolsep{1ex}}l}{\displaystyle  z_k(x)}\\
		\mathrm{ s.t.} &\displaystyle z_i(x) &\leq&  \varepsilon_i \quad \forall i = 1,\dots,p, i\neq k  \\
		&x &\in& X.
	\end{array}
\end{equation}

Every optimal solution of \eqref{eq:econs} is weakly efficient. If additionally the optimal solution of \eqref{eq:econs} is unique, it is an efficient solution of the multi-objective problem. Furthermore, all efficient solutions are optimal solutions of \eqref{eq:econs} for some vector \(\varepsilon\in\R^p\), i.e., all efficient solutions can be determined using the $\varepsilon$-constraint scalarization \citep[see, e.g.,][]{ehrgott05multicriteria}. We apply the $\varepsilon$-constraint method proposed in \citet{Kirlik2014a} that guarantees the efficiency of the generated solutions by using a two-stage approach.

We adaptively solve $\varepsilon$-constraint scalarizations of the root node problem to integer optimality to obtain possibly unsupported non-dominated points. The used $\varepsilon$ is obtained by the local upper bound $lu$, which spans the local hypervolume gap or the hypervolume of a search zone box of the corresponding node. Since we consider integer programs with integer coefficients, $\varepsilon$ can be chosen as $\varepsilon\coloneqq (lu_1-1,\dots, lu_p-1)^\top$. In the best case the obtained solution maps to a non-dominated point which has not been found before and can thus be added to the incumbent list. However, there is no guarantee that this point is an unsupported non-dominated point. It is not even guaranteed that the scalarized problem is feasible. Nevertheless, if a new non-dominated point is found the upper bound set is improved, which improves the fathoming rate.

The $\varepsilon$-constraint scalarization is applied between Step 2 and Step 3 of Algorithm~\ref{alg:BB}.

\subsubsection{Using Simple Lower Bound Sets}\label{subsec:SLB}
The methods proposed so far considered only scalarizations of the root node problem, such that the obtained objective space information can be integrated into every node and the corresponding subproblem, respectively. However, it is also possible to use scalarizations in subproblems of the branch and bound. Then, the obtained information does not hold for every node, but for all child nodes in the corresponding branch. Obtained solutions are efficient for the subproblems but in general not efficient for the underlying problem, since they might be dominated by solutions in other branches. 

So, solving IP~scalarizations in subproblems can be very time consuming and the obtained information might even be useless. However, solving a weighted sum scalarization to integer optimality compensates in some situations the additional costs. Instead of computing the complete lower bound set which might have many extreme supported points and facets, we adaptively solve a single weighted sum scalarization to integer optimality. As already discussed the level set of scalarization in an optimal solution is a valid lower bound on the objective values of the feasible solutions of the subproblem. So we save the time of computing the lower bound and use only the level set, obtained by solving the weighted sum scalarization, as the lower bound set. Of course this bound is weaker in most of the covered region compared to the lower bound obtained by solving the linear relaxation. This can be partially compensated by adding the so-called \emph{extreme facets} of the parent node to the lower bound set. By extreme facets we denote the facets of a lower bound set that are parallel to the axes. Thus, the simple lower bound set consists of $p+1$ facets. Then, the obtained inequality holds for every child node of this branch and can therefore be added as a constraint to the corresponding subproblems. Note that the test for redundancy of those inequalities is done in the same way as described in Section~\ref{subsub:warmstart}. In the best case, the obtained solution is efficient and has not been found before. Then, it can be added to the incumbent list and the upper bound is improved. 

A crucial component in this approach is the choice of the weight $\lambda \in \R^p$, which should be selected depending on the properties of the corresponding active node. Therefore,  we take into consideration all local upper bounds that were still located above the lower bound set in the parent node. Then, for each objective $k$, we choose among the considered local upper bounds the one with minimal objective value $lu_k$. Those $p$ points define a hyperplane $\mathcal{H}$. The dichotomic search approach uses this normal vector \(v\in\R^p\) of \(H\) as the weight $\lambda$. Unfortunately, for $p\geq 3$ this normal vector $v$ may have negative components and can not be used as weighting vector. Nevertheless, if \(v\geq 0\) is componentwise non-negative, we use it as weighting vector $\lambda=v$. Otherwise, we use the weight $\lambda=(1,\dots,1)^\top \in \R^p$.

\subsection{Algorithmic Control of the Presented Approaches}
In the previous sections improved components of multi-objective branch and bound algorithms are proposed. However, their algorithmic control is not covered yet, in particular, it is not specified when and how often IP~scalarizations should be solved. Since the costs of solving a scalarization to integer optimality are relatively high, this seems to be an important decision. In the first place, we want to obtain as much objective space information as possible. As a result, the lower and upper bound will be improved significantly. Those improved bounds increase the probability of fathoming by dominance and reduce the size of the search region. Both aspects will reduce the number of explored nodes, which will in turn reduce the total computational time. However, solving an excessive number of IP~scalarizations will increase the computation time. In addition, if IP~scalarizations are applied too often, there is an increasing chance that solved scalarizations do not provide new objective space information and are therefore redundant.

Obviously, there exists a trade-off between the decrease of the created nodes and the decrease of the total computation time. It is therefore necessary to find proper conditions, when to solve scalarizations to integer optimality. It is promising to already have good bounds in the early stages of the algorithm. This would lead to an improved fathoming rate from the beginning. Thus, we use a warmstarting approach, which solves IP~scalarizations in a preprocessing step. Like mentioned in Section~\ref{subsub:warmstart}, we use a predefined weight set for that. Preliminary numerical tests have shown, that already a small number of scalarizations solved to integer optimality have a high impact on the performance. We therefore use a predefined weight set $\Lambda$ with $|\Lambda | = p+1$. As weight vectors we use the standard unit vectors and additionally a weight vector with equal weights, i.e., $\lambda = (1,\dots,1)^\top \in \R^p$. Note that we add a small augmentation term to the vectors of the canonical basis to guarantee efficiency. 

Similar observations can be made, when the $\varepsilon$-constraint method is used to improve the upper bound. Preliminary tests have shown, that it is more promising to solve these IP~scalarizations in the early stages of the algorithm. This increases chances of obtaining unsupported non-dominated points early, which improves the upper bound set. Obviously, when the $\varepsilon$-constraint scalarization is applied too often, we probably waste a lot of time by solving multiple problems to integer optimality without a benefit. 

By using the simple lower bound approach instead of computing the complete lower bound set, some information is lost. The simple lower bound is in general weaker and there is no information about extreme points of the lower bound set. In the numerical tests, presented in Section~\ref{sec:res}, we will use the \emph{most-often fractional rule} (cf.\ Section~\ref{sec:relwork}). This rule can not be applied when the simple lower bound is used, as it is based on a single integer solution and thus there is no fractional variable. 
However, in this case we use the \emph{sum of ratios} branching rule that is presented in \citet{Bazgan2009solving} for multi-objective knapsack problems. They define ratio vectors $r_i \coloneqq (c^k_i \slash a_i)_{k=1,\dots,p}$ for all $i=1,\dots,n$ and sum the entries up, i.e., $sr_i\coloneqq \sum_{k=1}^p c^k_i \slash a_i$. Thereby, $a\in \R^n$ is the vector of constraint coefficents and $c^k$ is the $k$-th objective vector. The variable with the highest sum of ratios is branched on.
This can be considered as a direct extension of the basic single-objective branching rule for knapsack problems \citep[cf.][]{Kellerer2004knapsack} Note that this rule can also be applied to other problem classes like facility location problems or generalized assignment problems by interpreting the facility opening costs respectively the workload as weight. 

We decided to apply the simple lower bound approach at certain levels of the branch and bound tree.  This means, that we compute the simple lower bound instead of the complete bound set, when there is a certain number of fixed variables in the active node. Since there is in general a high number of nodes at deeper levels in the tree, a high amount of problems is solved to integer feasibility, resulting in possibly rising computation times. The efficiency and impact of the presented approaches are presented in the next section.

\section{Numerical Tests}\label{sec:res}
All presented algorithms were implemented in Julia 1.9.0 and the linear relaxations (for the lower bound set) were solved with Bensolve 2.1 \citep{Loehne2017the}. The scalarizations were solved to integer optimality with CPLEX 20.1. The numerical test runs were executed on a single core of a 3.20 GHz Intel\textsuperscript{\textregistered}  Core\texttrademark\ i7-8700 CPU  with 32~GB RAM. Note that the implementation of the proposed algorithms is publicly available \citep{Bauss2023gitimobb}.

We present different combinations of our presented approaches and compare their performance to a basic multi-objective branch and bound algorithm, which serves a baseline implementation. The results are evaluated regarding the average number of explored nodes and the average computational time over 10 instances. The time limit on solving a single instance is set to two hours.  

The basic branch and bound algorithm is constructed in the following way. 
\paragraph{Basic Branch and Bound}
\begin{itemize}
	\item \emph{Lower bound:} linear relaxation
	\item \emph{Upper bound:} incumbent list
	\item \emph{Node selection:} depth-first strategy
	\item \emph{Branching rule:} most-often fractional
\end{itemize}

Additionally to the basic branch and bound approach, we evaluate different combinations of the proposed approaches. These approaches show measurable impact in our test runs on different sets of problems. All considered branch and bound configurations are described in the following. 

\begin{itemize}
	\item \textbf{BB}. The basic branch and bound.
	\item \textbf{NS(.)}. The basic branch and bound, but with the dynamic node selection strategy, presented in Section~\ref{subsec:nodesel}. We distinguish between NS(LHG), when the local hypervolume gap is used, and NS(HSZ), when the hypervolume of the search zone box is considered.
	\item \textbf{WST}. Same procedure as NS(.), but using a warmstart of the bound sets.
	\item \textbf{EC}. Same procedure as WST. Additionally, $\varepsilon$-constraint scalarizations are solved. The scalarization is applied every $n$-th iteration within the first $p\,n^2$ iterations. 
	\item \textbf{SLB}. Same procedure as EC, but every fifth level of the branch and bound tree, the simple lower bound is considered, instead of solving the linear relaxation. 
	\item \textbf{+TE}. For every problem class, we consider the best performing approaches regarding the number of nodes respectively the total time. In those approaches if only 10 or less variables are free, we enumerate all $2^{10} = 1024$ solutions.
\end{itemize}

Note that the chosen parameters yield from preliminary numerical experiments on a different set of instances. Of course, they are not optimized and we do not change these parameter values for different problem classes, since we aim to show the impact of these approaches on a variety of problems. Finally, we present the considered problem classes and benchmark instances:
\begin{enumerate}[label={(\roman*)}]
\item Knapsack problems (KP) benchmark instances from \citet{Kirlik2014a}. The instances with 3 objectives and  40, 50, 60, 70 and 80 variables are solved. Additionally the instances with 4 objectives and 20, 30 and 40 variables are solved.
\item Uncapacitated facility location problems (UFLP) benchmark instances from \citet{Forget2022warm}. The instances with 3 objectives and 56, 72 and 90 variables are considered, as well as the instances with 4 objectives and 42 and 56 variables.
\item Capacitaded facility location problems (CFLP) instances from \citet{An2022a} and \citet{Bauss2023gitinstances}. We consider instances with 3 objectives and 65, 119 and 230 variables.
\item Generalized assignment problems (GAP) test instances from \citet{Bauss2023gitinstances}. The instances with 3 objectives and 48, 75 and 108 variables are solved. Additionally the instances with 4 objectives and 48 and 75 variables ware solved.
\end{enumerate}

\paragraph{Remarks Regarding the Implementation}
Due to numerical difficulties, we slightly adapt the implementation, as such it slightly differs from the presented approaches. The first change affects the inequalities, obtained by solving a weighted sum scalarization to integer optimality (cf.~Section~\ref{subsub:warmstart} and Section~\ref{subsec:SLB}). Originally, all constraint coefficients of the considered instances are integer. Nevertheless, the additional obtained inequalities, which are added to the corresponding subproblems, contain in general non-integer coefficients. Bensolve, which we use to compute our lower bound set, relies on the GLPK solver. Unfortunately, for harder and larger problems, the solver faces numerical issues and aborts the run of the algorithm in the worst case. To overcome this issues we round the constraint coefficients in the following way. 

Let $\bar{a}_1x_1+,\dots,+\bar{a}_nx_n \geq \bar{b}$, with $ \bar{a} \in \R^n$ and $\bar{b}\in\R$, be an inequality obtained during the algorithm. Then, the modified inequality is given by $\lceil\bar{a}_1\rceil x_1+,\dots,\lceil +\bar{a}_n\rceil x_n \geq \lfloor \bar{b}\rfloor$. Obviously, this constraint is weaker in general, but it is necessary to overcome the numerical issues. 

The second change also concerns the usage of the simple lower bounds (cf.~Section~\ref{subsec:SLB}). The main motivation of using this simple lower bounds was to save computation time by not using Bensolve to compute the complete lower bound set which may consist of a large number of facets. But although we use CPLEX to solve just a single scalarization to integer optimality instead of the complete lower bound set, this is unfortunately much slower in the majority of the considered instances. This is due to the rather slow interface to CPLEX in Julia. Especially the problem building takes a relatively large amount of time. To partially overcome this problem, we limit the CPLEX computation time to $\frac{1}{10}$-th of the time needed to compute the lower bound with Bensolve in the root node. If the program is not solved to optimality within the the given time, the best feasible solution found so far is treated like the optimal solution to possibly update the upper bound and the current best lower bound can be used as simple lower bound. Therefore it is also possible to use objective space information, although they might be worse. If no feasible solution is found within the time limit and infeasibility is not proven, no additional objective information can be added to the corresponding node. 

\begin{table}
	\scriptsize
	\hspace{-1.2cm}
	\subfloat[Knapsack problem with $n=40$ variables and $p=3$ objectives	\label{Knap3-40}]{
		\begin{filecontents*}{Results/K40.csv}
approach,nodes,time (s),IPs,solved
BB,138365.8,97.55,0.0,10
NS(LHG),49898.4,76.85,0.0,10
WST,48157.8,73.00,4.0,10
EC,45384.2,68.99,39.3,10
SLB,51869.0,100.17,9203.5,10
EC+TE,41948.0,68.76,38.4,10
SLB+TE,33993.0,106.21,6838.4,10
\end{filecontents*}
\pgfplotstabletypeset[
		col sep=comma,
		string type,
		every head row/.style={%
			before row={\hline
				\multicolumn{5}{|l|}{\rule{0pt}{1em}knapsack problem, $p=3,n=40$}  \\\hline
			},
			after row=\hline
		},
		every last row/.style={after row=\hline},
		columns/approach/.style={column name = \centering approach, column type=|l},
		columns/nodes/.style={column name= \centering nodes, column type=|R{1.2cm}},
		columns/time (s)/.style={column name= \centering time (s), column type=|R{1.2cm}},
		columns/IPs/.style={column name=\centering \#IPs, column type= |R{1cm}},
		columns/solved/.style={column name=  solved, column type=|c|},
		]{Results/K40.csv}
	}
	\hspace{0.25cm}
	\subfloat[Knapsack problem with $n=50$ variables and $p=3$ objectives	\label{Knap3-50}]{
		\begin{filecontents*}{Results/K50.csv}
approach,nodes,time (s),IPs,solved
BB,391170.2,398.55,0.0,10
NS(LHG),112975.4,317.70,0.0,10
WST,110409.0,299.95,4.0,10
EC,101027.0,254.48,50.6,10
SLB,97466.2,391.33,17479.4,10
EC+TE,98385.4,256.40,49.4,10
SLB+TE,94661.0,375.52,17492.8,10
\end{filecontents*}
\pgfplotstabletypeset[
		col sep=comma,
		string type,
		every head row/.style={%
			before row={\hline
				\multicolumn{5}{|l|}{\rule{0pt}{1em}knapsack problem, $p=3,n=50$}  \\\hline
			},
			after row=\hline
		},
		every last row/.style={after row=\hline},
		columns/approach/.style={column name = \centering approach, column type=|l},
		columns/nodes/.style={column name= \centering nodes, column type=|R{1.2cm}},
		columns/time (s)/.style={column name= \centering time (s), column type=|R{1.2cm}},
		columns/IPs/.style={column name=\centering \#IPs, column type= |R{1cm}},
		columns/solved/.style={column name=  solved, column type=|c|},
		]{Results/K50.csv}
	}

	\hspace{-1.2cm}
	\subfloat[Knapsack problem with $n=60$ variables and $p=3$ objectives	\label{Knap3-60}]{
		\begin{filecontents*}{Results/K60.csv}
approach,nodes,time (s),IPs,solved
BB,1201949.8,1794.81,0.0,10
NS(LHG),321465.3,2591.44,0.0,8
WST,320631.0,2550.18,4.0,9
EC,313461.2,2103.35,62.1,10
SLB,247905.9,1946.32,45446.9,10
EC+TE,279708.3,1496.39,63.1,10
SLB+TE,276558.0,2082.41,49992.1,10
\end{filecontents*}
\pgfplotstabletypeset[
		col sep=comma,
		string type,
		every head row/.style={%
			before row={\hline
				\multicolumn{5}{|l|}{\rule{0pt}{1em}knapsack problem, $p=3,n=60$}  \\\hline
			},
			after row=\hline
		},
		every last row/.style={after row=\hline},
		columns/approach/.style={column name = \centering approach, column type=|l},
		columns/nodes/.style={column name= \centering nodes, column type=|R{1.2cm}},
		columns/time (s)/.style={column name= \centering time (s), column type=|R{1.2cm}},
		columns/IPs/.style={column name=\centering \#IPs, column type= |R{1cm}},
		columns/solved/.style={column name=  solved, column type=|c|},
		]{Results/K60.csv}
	}
	\hspace{0.25cm}
	\subfloat[Knapsack problem with $n=70$ variables and $p=3$ objectives	\label{Knap3-70}]{
		\begin{filecontents*}{Results/K70.csv}
approach,nodes,time (s),IPs,solved
BB,2470128.8,4312.96,0.0,7
NS(LHG),447631.9,4533.68,0.0,7
WST,447604.2,4489.07,4.0,7
EC,432262.9,4156.09,67.2,7
SLB,440161.9,4532.67,79311.0,7
EC+TE,439530.9,3828.58,67.8,8
SLB+TE,436677.6,4572.73,79311.2,8
\end{filecontents*}
\pgfplotstabletypeset[
		col sep=comma,
		string type,
		every head row/.style={%
			before row={\hline
				\multicolumn{5}{|l|}{\rule{0pt}{1em}knapsack problem, $p=3,n=70$}  \\\hline
			},
			after row=\hline
		},
		every last row/.style={after row=\hline},
		columns/approach/.style={column name = \centering approach, column type=|l},
		columns/nodes/.style={column name= \centering nodes, column type=|R{1.2cm}},
		columns/time (s)/.style={column name= \centering time (s), column type=|R{1.2cm}},
		columns/IPs/.style={column name=\centering \#IPs, column type= |R{1cm}},
		columns/solved/.style={column name=  solved, column type=|c|},
		]{Results/K70.csv}
	}
	
	\hspace{-1.2cm}
	\subfloat[Knapsack problem with $n=80$ variables and $p=3$ objectives	\label{Knap3-80}]{
		\begin{filecontents*}{Results/K80.csv}
approach,nodes,time (s),IPs,solved
BB,3433775.1,7200.00,0.0,0
NS(LHG),376535.2,7200.00,0.0,0
WS,392186.7,7200.00,4.0,0
EC,429098.6,7200.00,88.8,0
SLB,500545.7,7158.29,87007.5,1
EC+TE,418424.2,7200.00,89.1,0
SLB+TE,505265.6,7105.14,88435.9,1
\end{filecontents*}
\pgfplotstabletypeset[
		col sep=comma,
		string type,
		every head row/.style={%
			before row={\hline
				\multicolumn{5}{|l|}{\rule{0pt}{1em}knapsack problem, $p=3,n=80$}  \\\hline
			},
			after row=\hline
		},
		every last row/.style={after row=\hline},
		columns/approach/.style={column name = \centering approach, column type=|l},
		columns/nodes/.style={column name= \centering nodes, column type=|R{1.2cm}},
		columns/time (s)/.style={column name= \centering time (s), column type=|R{1.2cm}},
		columns/IPs/.style={column name=\centering \#IPs, column type= |R{1cm}},
		columns/solved/.style={column name=  solved, column type=|c|},
		]{Results/K80.csv}
	}
	\hspace{0.25cm}
	\subfloat[Knapsack problem with $n=20$ variables and $p=4$ objectives	\label{Knap4-20}]{
		\begin{filecontents*}{Results/K4-20.csv}
approach,nodes,time (s),IPs,solved
BB,9073.6,7.76,0.0,10
NS(LHG),4398.0,6.19,0.0,10
WST,4230.8,6.33,5.0,10
EC,4185.6,7.12,38.7,10
SLB,4289.8,15.57,752.5,10
EC+TE,899.0,4.90,27.7,10
SLB+TE,742.8,8.73,303.4,10
\end{filecontents*}
\pgfplotstabletypeset[
		col sep=comma,
		string type,
		every head row/.style={%
			before row={\hline
				\multicolumn{5}{|l|}{\rule{0pt}{1em}knapsack problem, $p=4,n=20$}  \\\hline
			},
			after row=\hline
		},
		every last row/.style={after row=\hline},
		columns/approach/.style={column name = \centering approach, column type=|l},
		columns/nodes/.style={column name= \centering nodes, column type=|R{1.2cm}},
		columns/time (s)/.style={column name= \centering time (s), column type=|R{1.2cm}},
		columns/IPs/.style={column name=\centering \#IPs, column type= |R{1cm}},
		columns/solved/.style={column name=  solved, column type=|c|},
		]{Results/K4-20.csv}
	}
	
	\hspace{-1.2cm}
	\subfloat[Knapsack problem with $n=30$ variables and $p=4$ objectives	\label{Knap4-30}]{
		\begin{filecontents*}{Results/K4-30.csv}
approach,nodes,time (s),IPs,solved
BB,53531.0,90.03,0.0,10
NS(LHG),21621.0,78.95,0.0,10
WST,19798.3,68.42,5.0,10
EC,19107.2,67.78,59.2,10
SLB,26117.0,206.31,4699.9,10
EC+TE,15618.0,73.85,59.3,10
SLB+TE,14554.0,135.63,3428.7,10
\end{filecontents*}
\pgfplotstabletypeset[
		col sep=comma,
		string type,
		every head row/.style={%
			before row={\hline
				\multicolumn{5}{|l|}{\rule{0pt}{1em}knapsack problem, $p=4,n=30$}  \\\hline
			},
			after row=\hline
		},
		every last row/.style={after row=\hline},
		columns/approach/.style={column name = \centering approach, column type=|l},
		columns/nodes/.style={column name= \centering nodes, column type=|R{1.2cm}},
		columns/time (s)/.style={column name= \centering time (s), column type=|R{1.2cm}},
		columns/IPs/.style={column name=\centering \#IPs, column type= |R{1cm}},
		columns/solved/.style={column name=  solved, column type=|c|},
		]{Results/K4-30.csv}
	}
	\hspace{0.25cm}
	\subfloat[Knapsack problem with $n=40$ variables and $p=4$ objectives	\label{Knap4-40}]{
		\begin{filecontents*}{Results/K4-40.csv}
approach,nodes,time (s),IPs,solved
BB,325416.1,2070.93,0.0,8
NS(LHG),115926.3,1815.93,0.0,8
WST,114715.6,1562.56,5.0,8
EC,110539.0,1511.55,72.4,8
SLB,85321.0,1410.12,14958.6,8
EC+TE,111540.9,2235.66,68.6,8
SLB+TE,73076.1,2974.74,11223.3,8
\end{filecontents*}
\pgfplotstabletypeset[
		col sep=comma,
		string type,
		every head row/.style={%
			before row={\hline
				\multicolumn{5}{|l|}{\rule{0pt}{1em}knapsack problem, $p=4,n=40$}  \\\hline
			},
			after row=\hline
		},
		every last row/.style={after row=\hline},
		columns/approach/.style={column name = \centering approach, column type=|l},
		columns/nodes/.style={column name= \centering nodes, column type=|R{1.2cm}},
		columns/time (s)/.style={column name= \centering time (s), column type=|R{1.2cm}},
		columns/IPs/.style={column name=\centering \#IPs, column type= |R{1cm}},
		columns/solved/.style={column name=  solved, column type=|c|},
		]{Results/K4-40.csv}
	}
	
	\caption{Numerical results on multi-objective knapsack instances of \citet{Kirlik2014a}.}
\end{table}		

\paragraph{Results on Knapsack Problems}
The numerical results show that the dynamic node selection strategy based on the local hypervolume gap has a large impact on the average number of explored nodes. Note that preliminary tests on a different set of knapsack problems have shown, that this dynamic strategy works better than the strategy based on the search zone box. With the chosen node selection strategy we can reduce the number of explored nodes by up to $71.1\%$ (Table~\ref{Knap3-50}), in problem sizes where all $10$ instances are solved, and up to $81.9\%$ (Table~\ref{Knap3-70}) in instance sizes where the same amount of problems is solved. Since the computation of the local hypervolume gap can be expensive, the total computation time can only be reduced by up to $21.6\%$ (Table~\ref{Knap3-40}). Note, that there are also cases, where the total computation time increases, although the number of explored nodes decreases (Table~\ref{Knap3-60} and \ref{Knap3-70}). It is not surprising, that this occurs in the instances with a larger amount of variables, since there are possibly more non-dominated points and corresponding local upper bounds, which need to be considered during the gap computation. 

The approaches WST and EC reduce the average number of explored nodes in the majority of the considered instances. Nevertheless, the impact on the performance is more significant when EC is used. The usage of SLB can, especially for instances with a larger amount of variables, reduce the computation time and the number of explored nodes. In Table~\ref{Knap3-80}, it is shown, that using SLB allows us to solve more instances in the given time limit, which is an improvement compared to the other approaches. Nevertheless, especially for smaller instance sizes, using SLB increases the computational time a lot. This is due to the fact, that the CPLEX interface is rather slow in Julia.

The best working approaches for the considered benchmark instances are EC+TE and SLB+TE. Regarding the average number of explored nodes SLB+TE seems to be the best choice, since this approach creates the least amount of nodes in the majority of the experiments. It is possible to reduce the number of nodes by up to $91.8\%$ (Table~\ref{Knap4-20}), respectively $79.4\%$ (Table~\ref{Knap3-60}) if we omit the instance size $p=4,n=20$, since it highly benefits from the enumeration. The best choice for KP regarding the total computation time seems to be EC+TE, since it is the fastest approach in the majority of the solved instances. The runtime can be reduced by up to $36.9\%$ (Table~\ref{Knap4-20}), respectively $36.2\%$ (Table~\ref{Knap3-50}).

\begin{table}
	\scriptsize
	\hspace{-1.2cm}
	\subfloat[Uncapacitated facility location problem with $n=56$ variables and $p=3$ objectives	\label{UFLP3-7}]{
		\begin{filecontents*}{Results/FL7.csv}
approach,nodes,time (s),IPs,solved
BB,152349.4,216.07,0.0,10
NS(HSZ),100222.6,217.16,0.0,10
WST,97243.4,204.34,4.0,10
EC,97649.4,205.35,39.0,10
SLB,95000.8,258.92,17888.5,10
WST+TE,96318.8,216.03,4.0,10
SLB+TE,95795.9,318.66,18159.4,10
\end{filecontents*}
\pgfplotstabletypeset[
		col sep=comma,
		string type,
		every head row/.style={%
			before row={\hline
				\multicolumn{5}{|l|}{\rule{0pt}{1em}Uncapacitated facility location problem, $p=3,n=56$}  \\\hline
			},
			after row=\hline
		},
		every last row/.style={after row=\hline},
		columns/approach/.style={column name = \centering approach, column type=|l},
		columns/nodes/.style={column name= \centering nodes, column type=|R{1.2cm}},
		columns/time (s)/.style={column name= \centering time (s), column type=|R{1.2cm}},
		columns/IPs/.style={column name=\centering \#IPs, column type= |R{1cm}},
		columns/solved/.style={column name=  solved, column type=|c|},
		]{Results/FL7.csv}
	}
\hspace{0.25cm}
\subfloat[Uncapacitated facility location problem with $n=72$ variables and $p=3$ objectives	\label{UFLP3-72}]{
	\begin{filecontents*}{Results/FL8.csv}
approach,nodes,time (s),IPs,solved
BB,482483.6,1321.29,0.0,10
NS(HSZ),277068.6,1249.20,0.0,10
WST,310494.8,1377.52,4.0,10
EC,310802.8,1378.46,42.1,10
SLB,295379.0,1694.93,54184.4,10
WST+TE,310357.8,1310.23,4.0,10
SLB+TE,295245.3,1729.24,54187.9,10
\end{filecontents*}
\pgfplotstabletypeset[
	col sep=comma,
	string type,
	every head row/.style={%
		before row={\hline
			\multicolumn{5}{|l|}{\rule{0pt}{1em}Uncapacitated facility location problem, $p=3,n=72$}  \\\hline
		},
		after row=\hline
	},
	every last row/.style={after row=\hline},
	columns/approach/.style={column name = \centering approach, column type=|l},
	columns/nodes/.style={column name= \centering nodes, column type=|R{1.2cm}},
	columns/time (s)/.style={column name= \centering time (s), column type=|R{1.2cm}},
	columns/IPs/.style={column name=\centering \#IPs, column type= |R{1cm}},
	columns/solved/.style={column name=  solved, column type=|c|},
	]{Results/FL8.csv}
}

\hspace{-1.2cm}
\subfloat[Uncapacitated facility location problem with $n=90$ variables and $p=3$ objectives	\label{UFLP3-90}]{
	\begin{filecontents*}{Results/FL9.csv}
approach,nodes,time (s),IPs,solved
BB,1604980.2,7144.28,0.0,2
NS(HSZ),757425.0,6742.37,0.0,7
WST,712808.9,6328.50,4.0,8
EC,733769.6,6310.58,48.1,8
SLB,715377.1,6814.74,105225.6,6
WST+TE,731303.1,6326.04,4.0,8
SLB+TE,731190.3,6623.71,107625.1,7
\end{filecontents*}
\pgfplotstabletypeset[
	col sep=comma,
	string type,
	every head row/.style={%
		before row={\hline
			\multicolumn{5}{|l|}{\rule{0pt}{1em}Uncapacitated facility location problem, $p=3,n=90$}  \\\hline
		},
		after row=\hline
	},
	every last row/.style={after row=\hline},
	columns/approach/.style={column name = \centering approach, column type=|l},
	columns/nodes/.style={column name= \centering nodes, column type=|R{1.2cm}},
	columns/time (s)/.style={column name= \centering time (s), column type=|R{1.2cm}},
	columns/IPs/.style={column name=\centering \#IPs, column type= |R{1cm}},
	columns/solved/.style={column name=  solved, column type=|c|},
	]{Results/FL9.csv}
}
	\hspace{0.25cm}
\subfloat[Uncapacitated facility location problem with $n=42$ variables and $p=4$ objectives	\label{UFLP4-42}]{
	\begin{filecontents*}{Results/FL64.csv}
approach,nodes,time (s),IPs,solved
BB,61048.4,293.36,0.0,10
NS(HSZ),48404.8,272.20,0.0,10
WST,40328.0,226.39,5.0,10
EC,40947.6,228.46,27.4,10
SLB,42194.3,434.08,7406.1,10
WST+TE,37312.8,273.25,5.0,10
SLB+TE,37446.8,432.76,6918.8,10
\end{filecontents*}
\pgfplotstabletypeset[
	col sep=comma,
	string type,
	every head row/.style={%
		before row={\hline
			\multicolumn{5}{|l|}{\rule{0pt}{1em}Uncapacitated facility location problem, $p=4,n=42$}  \\\hline
		},
		after row=\hline
	},
	every last row/.style={after row=\hline},
	columns/approach/.style={column name = \centering approach, column type=|l},
	columns/nodes/.style={column name= \centering nodes, column type=|R{1.2cm}},
	columns/time (s)/.style={column name= \centering time (s), column type=|R{1.2cm}},
	columns/IPs/.style={column name=\centering \#IPs, column type= |R{1cm}},
	columns/solved/.style={column name=  solved, column type=|c|},
	]{Results/FL64.csv}
}

\hspace{-1.2cm}
\subfloat[Uncapacitated facility location problem with $n=56$ variables and $p=4$ objectives	\label{UFLP4-56}]{
	\begin{filecontents*}{Results/FL74.csv}
approach,nodes,time (s),IPs,solved
BB,435160.6,6175.10,0.0,5
NS(HSZ),313120.4,5745.16,0.0,8
WST,264705.2,4309.92,5.0,10
EC,266740.2,4324.72,34.6,10
SLB,180643.0,5296.01,29414.9,6
WST+TE,262173.8,4721.79,5.0,10
SLB+TE,165112.8,5484.16,30455.25,6
\end{filecontents*}
\pgfplotstabletypeset[
	col sep=comma,
	string type,
	every head row/.style={%
		before row={\hline
			\multicolumn{5}{|l|}{\rule{0pt}{1em}Uncapacitated facility location problem, $p=4,n=56$}  \\\hline
		},
		after row=\hline
	},
	every last row/.style={after row=\hline},
	columns/approach/.style={column name = \centering approach, column type=|l},
	columns/nodes/.style={column name= \centering nodes, column type=|R{1.2cm}},
	columns/time (s)/.style={column name= \centering time (s), column type=|R{1.2cm}},
	columns/IPs/.style={column name=\centering \#IPs, column type= |R{1cm}},
	columns/solved/.style={column name=  solved, column type=|c|},
	]{Results/FL74.csv}
}
	\caption{Numerical results for multi-objective uncapacitated facility location instances of \citet{Forget2022warm}. }
\end{table}		

\paragraph{Results on Uncapacitated Facility Location Problems}
Similar to the results of KP, the dynamic node selection strategy has a significant impact on the number of explored nodes and the total computation time. In contrast to the knapsack problems, preliminary tests on a different set of UFLP have shown, that the dynamic node selection strategy based on the search zone box works better. The reason for this might be the larger amount of non-dominated points, that cause an even larger amount of local upper bounds. Again, a high amount of local upper bounds increases the time needed to compute the gap between lower and upper bound set. Therefore, the gap measure based on the search zone box performs better, since its computation is way faster. The average number of explored nodes can be reduced by up to $42.2\%$ (Table~\ref{UFLP3-72}) and the total computation time can be reduced by up to $7.2\%$ (Table~\ref{UFLP3-72}). Furthermore, disregarding the improvements in explored nodes and runtime, there are also improvements regarding the number of solved instances within the given time limit of two hours. For the instance size $n=90, p=3$ the amount of solved instances is improved from $2$ to $7$ (Table~\ref{UFLP3-90}) and for the instance size $n=56,p=4$ the amount of solved instances is improved from $5$ to $8$ (Table~\ref{UFLP4-56}).

The WST method reduces the number of explored nodes and the total computation time, in comparison with NS(HSZ), in nearly all considered instance sizes. Table~\ref{UFLP3-72} is the only exception. Although EC yields similar results like WST, it is slightly worse in the majority of the benchmark instances. Therefore, there is no need to solve the additional $\varepsilon$-constraint scalarizations, since it is outperformed by WST. Although it seems that the number of explored nodes can be further reduced with SLB, the computation time increases. This even results in a reduced number of solved instances within the time limit (Table~\ref{UFLP3-90} and Table~\ref{UFLP4-56}).

Regarding the total computation time WST is the best approach to solve uncapacitated facility problems, since it is the fastest in the majority of the instances. In the best case, the computation time is reduced by up to $30.2 \%$ (Table~\ref{UFLP4-56}), where at the same time the number of solved instances is increased from $5$ to $10$. Regarding the number of explored nodes WST+TE seems to be a good choice. In the majority of the benchmark instances, this approach yields the lowest number of explored nodes. The average number of explored nodes can be reduced by up to $55.5\%$ in the best case (Table~\ref{UFLP3-90}).


\begin{table}
	\scriptsize
	\hspace{-1.2cm}
	\subfloat[Capacitated facility location problem with $n=65$ variables and $p=3$ objectives	\label{CFLP3-5}]{
		\begin{filecontents*}{Results/CFL5-10.csv}
approach,nodes,time (s),IPs,solved
BB,59890.6,48.40,0.0,10
NS(HSZ),49007.6,46.31,0.0,10
WST,48857.0,47.18,4.0,10
EC,48839.6,45.89,20.6,10
SLB,48587.2,55.70,9195.8,10
EC+TE,48839.6,49.50,20.6,10
SLB+TE,47988.4,55.36,9203.3,10
\end{filecontents*}
\pgfplotstabletypeset[
		col sep=comma,
		string type,
		every head row/.style={%
			before row={\hline
				\multicolumn{5}{|l|}{\rule{0pt}{1em}Capacitated facility location problem, $p=3,n=65$}  \\\hline
			},
			after row=\hline
		},
		every last row/.style={after row=\hline},
		columns/approach/.style={column name = \centering approach, column type=|l},
		columns/nodes/.style={column name= \centering nodes, column type=|R{1.2cm}},
		columns/time (s)/.style={column name= \centering time (s), column type=|R{1.2cm}},
		columns/IPs/.style={column name=\centering \#IPs, column type= |R{1cm}},
		columns/solved/.style={column name=  solved, column type=|c|},
		]{Results/CFL5-10.csv}
	}
	\hspace{0.25cm}
	\subfloat[Capacitated facility location problem with $n=119$ variables and $p=3$ objectives	\label{CFLP3-7}]{
		\begin{filecontents*}{Results/CFL7-14.csv}
approach,nodes,time (s),IPs,solved
BB,460374.4,1134.55,0.0,10
NS(HSZ),321311.6,919.95,0.0,10
WST,327431.8,939.44,4.0,10
EC,326869.2,919.65,23.8,10
SLB,345477.6,1009.03,62699.5,10
EC+TE,326869.2,942.93,23.8,10
SLB+TE,345477.6,1012.17,62699.5,10
\end{filecontents*}
\pgfplotstabletypeset[
		col sep=comma,
		string type,
		every head row/.style={%
			before row={\hline
				\multicolumn{5}{|l|}{\rule{0pt}{1em}Capacitated facility location problem, $p=3,n=119$}  \\\hline
			},
			after row=\hline
		},
		every last row/.style={after row=\hline},
		columns/approach/.style={column name = \centering approach, column type=|l},
		columns/nodes/.style={column name= \centering nodes, column type=|R{1.2cm}},
		columns/time (s)/.style={column name= \centering time (s), column type=|R{1.2cm}},
		columns/IPs/.style={column name=\centering \#IPs, column type= |R{1cm}},
		columns/solved/.style={column name=  solved, column type=|c|},
		]{Results/CFL7-14.csv}
	}
	
	\hspace{-1.2cm}
	\subfloat[Capacitated facility location problem with $n=230$ variables and $p=3$ objectives	\label{CFLP3-10}]{
		\begin{filecontents*}{Results/CFL10-20.csv}
approach,nodes,time (s),IPs,solved
BB,660929.5,7200.00,0.0,0
NS(HSZ),358055.8,7200.00,0.0,0
WST,362319.4,7200.00,0.0,0
EC,387534.1,7200.00,16.8,0
SLB,453839.3,7200.00,80912.8,0
EC+TE,391868.6,7200.00,16.8,0
SLB+TE,427363.5,7200.00,83468.0,0
\end{filecontents*}
\pgfplotstabletypeset[
		col sep=comma,
		string type,
		every head row/.style={%
			before row={\hline
				\multicolumn{5}{|l|}{\rule{0pt}{1em}Capacitated facility location problem, $p=3,n=230$}  \\\hline
			},
			after row=\hline
		},
		every last row/.style={after row=\hline},
		columns/approach/.style={column name = \centering approach, column type=|l},
		columns/nodes/.style={column name= \centering nodes, column type=|R{1.2cm}},
		columns/time (s)/.style={column name= \centering time (s), column type=|R{1.2cm}},
		columns/IPs/.style={column name=\centering \#IPs, column type= |R{1cm}},
		columns/solved/.style={column name=  solved, column type=|c|},
		]{Results/CFL10-20.csv}
	}

	\caption{Numerical results for multi-objective capacitated facility location instances of \citet{An2022a} and \citet{Bauss2023gitinstances}.}
\end{table}

\paragraph{Results of Capacitated Facility Location Problems}
Similar to the uncapacitated facility location problem, preliminary tests on a different set of instances have shown that the dynamic search strategy based on the search zone boxes are computationally more efficient. The number of nodes and the computation time are improved with NS(HSZ). By using the suggested node selection strategy, the number of explored nodes can be reduced by up $30.2\%$ (Table~\ref{CFLP3-7}) and the total computation time improves by up to $18.9\%$ (Table~\ref{CFLP3-7}). Using the warmstarting of bound sets (WST), has nearly no impact on the number of explored nodes but the computational time increases. By using EC, the number of nodes is reduced marginally, but there is an improvement regarding the computation time, compared to WST.  By also using the simple lower bounds, i.e., using the SLB approach, the number of nodes and the computation time increases. Nevertheless, compared to our baseline implementation (BB) there are improvements regarding the number of nodes and there can also be an improvement regarding the runtime.

The best approach regarding the computational time seems to be EC. In the best case, the time can be reduced by $18.9\%$ (Table~\ref{CFLP3-7}). Regarding the number of explored nodes, there seems to be no consistency caused by the small amount of benchmark instances. NS(HSZ), SLB+TE and EC seem to be good choices to decrease the number of explored nodes, where those can be reduced by $30.2\%$ in the best case.


\begin{table}
	\scriptsize
	\hspace{-1.2cm}
	\subfloat[Generalized assignment problem with $n=48$ variables and $p=3$ objectives	\label{GAP3-12}]{
		\begin{filecontents*}{Results/GAP12.csv}
approach,nodes,time (s),IPs,solved
BB,27890.0,16.76,0.0,10
NS(HSZ),20608.0,14.57,0.0,10
WST,18458.4,13.17,4.0,10
EC,18506.4,13.96,32.7,10
SLB,18567.8,20.51,3516.7,10
WST+TE,17469.0,14.08,4.0,10
EC+TE,17501.0,14.91,34.3,10
SLB+TE,17895.8,21.27,3336.5,10
\end{filecontents*}
\pgfplotstabletypeset[
		col sep=comma,
		string type,
		every head row/.style={%
			before row={\hline
				\multicolumn{5}{|l|}{\rule{0pt}{1em}Generalized assignment problem, $p=3,n=48$}  \\\hline
			},
			after row=\hline
		},
		every last row/.style={after row=\hline},
		columns/approach/.style={column name = \centering approach, column type=|l},
		columns/nodes/.style={column name= \centering nodes, column type=|R{1.2cm}},
		columns/time (s)/.style={column name= \centering time (s), column type=|R{1.2cm}},
		columns/IPs/.style={column name=\centering \#IPs, column type= |R{1cm}},
		columns/solved/.style={column name=  solved, column type=|c|},
		]{Results/GAP12.csv}
	}
	\hspace{0.25cm}
	\subfloat[Generalized assignment problem with $n=75$ variables and $p=3$ objectives	\label{GAP3-15}]{
		\begin{filecontents*}{Results/GAP15.csv}
approach,nodes,time (s),IPs,solved
BB,170359.4,220.20,0.0,10
NS(HSZ),98703.4,162.00,0.0,10
WST,95179.0,154.56,4.0,10
EC,95760.4,157.35,44.5,10
SLB,66464.0,119.34,12716.6,10
WST+TE,94572.0,166.30,4.0,10
EC+TE,95134.0,168.09,45.1,10
SLB+TE,61568.0,113.74,11676.6,10
\end{filecontents*}
\pgfplotstabletypeset[
		col sep=comma,
		string type,
		every head row/.style={%
			before row={\hline
				\multicolumn{5}{|l|}{\rule{0pt}{1em}Generalized assignment problem, $p=3,n=75$}  \\\hline
			},
			after row=\hline
		},
		every last row/.style={after row=\hline},
		columns/approach/.style={column name = \centering approach, column type=|l},
		columns/nodes/.style={column name= \centering nodes, column type=|R{1.2cm}},
		columns/time (s)/.style={column name= \centering time (s), column type=|R{1.2cm}},
		columns/IPs/.style={column name=\centering \#IPs, column type= |R{1cm}},
		columns/solved/.style={column name=  solved, column type=|c|},
		]{Results/GAP15.csv}
	}
	
	\hspace{-1.2cm}
	\subfloat[Generalized assignment problem with $n=108$ variables and $p=3$ objectives	\label{GAP3-18}]{
		\begin{filecontents*}{Results/GAP18.csv}
approach,nodes,time (s),IPs,solved
BB,1041528.2,2566.35,0.0,10
NS(HSZ),484977.8,1569.04,0.0,10
WST,467268.8,1527.91,4.0,10
EC,469200.1,1534.35,53.7,10
SLB,387682.2,1603.55,74288.6,10
WST+TE,464769.0,1632.80,4.0,10
EC+TE,466327.2,1650.48,55.3,10
SLB+TE,367327.0,1746.80,68916.3,10
\end{filecontents*}
\pgfplotstabletypeset[
		col sep=comma,
		string type,
		every head row/.style={%
			before row={\hline
				\multicolumn{5}{|l|}{\rule{0pt}{1em}Generalized assignment problem, $p=3,n=108$}  \\\hline
			},
			after row=\hline
		},
		every last row/.style={after row=\hline},
		columns/approach/.style={column name = \centering approach, column type=|l},
		columns/nodes/.style={column name= \centering nodes, column type=|R{1.2cm}},
		columns/time (s)/.style={column name= \centering time (s), column type=|R{1.2cm}},
		columns/IPs/.style={column name=\centering \#IPs, column type= |R{1cm}},
		columns/solved/.style={column name=  solved, column type=|c|},
		]{Results/GAP18.csv}
	}
	\hspace{0.25cm}
	\subfloat[Generalized assignment problem with $n=48$ variables and $p=4$ objectives	\label{GAP4-12}]{
		\begin{filecontents*}{Results/GAP12-4.csv}
approach,nodes,time,IPs,solved
BB,145278.4,529.94,0.0,10
NS(HSZ),99512.2,494.75,0.0,10
WST,85777.2,436.29,5.0,10
EC,77102.8,325.46,26.6,10
SLB,74366.3,386.94,13771.3,10
WST+TE,78313.0,494.43,5.0,10
EC+TE,78622.2,437.83,29.0,10
SLB+TE,73206.1,797.21,13150.0,10
\end{filecontents*}
\pgfplotstabletypeset[
		col sep=comma,
		string type,
		every head row/.style={%
			before row={\hline
				\multicolumn{5}{|l|}{\rule{0pt}{1em}Generalized assignment problem, $p=4,n=48$}  \\\hline
			},
			after row=\hline
		},
		every last row/.style={after row=\hline},
		columns/approach/.style={column name = \centering approach, column type=|l},
		columns/nodes/.style={column name= \centering nodes, column type=|R{1.2cm}},
		columns/time (s)/.style={column name= \centering time (s), column type=|R{1.2cm}},
		columns/IPs/.style={column name=\centering \#IPs, column type= |R{1cm}},
		columns/solved/.style={column name=  solved, column type=|c|},
		]{Results/GAP12-4.csv}
	}
	
	\hspace{-1.2cm}
	\subfloat[Generalized assignment problem with $n=75$ variables and $p=4$ objectives	\label{GAP4-15}]{
		\begin{filecontents*}{Results/GAP15-4.csv}
approach,nodes,time (s),IPs,solved
BB,836528.8,6292.90,0.0,4
NS(HSZ),490997.9,6185.76,0.0,5
WST,466591.2,6128.25,5.0,5
EC,466653.3,5984.18,33.4,5
SLB,297399.0,6423.18,54521.0,5
WST+TE,463480.6,6013.66,5.0,6
EC+TE,466678.0,6004.96,35.0,6
SLB+TE,282099.1,6925.44,55286.5,5
\end{filecontents*}
\pgfplotstabletypeset[
		col sep=comma,
		string type,
		every head row/.style={%
			before row={\hline
				\multicolumn{5}{|l|}{\rule{0pt}{1em}Generalized assignment problem, $p=4,n=75$}  \\\hline
			},
			after row=\hline
		},
		every last row/.style={after row=\hline},
		columns/approach/.style={column name = \centering approach, column type=|l},
		columns/nodes/.style={column name= \centering nodes, column type=|R{1.2cm}},
		columns/time (s)/.style={column name= \centering time (s), column type=|R{1.2cm}},
		columns/IPs/.style={column name=\centering \#IPs, column type= |R{1cm}},
		columns/solved/.style={column name=  solved, column type=|c|},
		]{Results/GAP15-4.csv}
	}
	\caption{Numerical results for multi-objective generalized assignment instances of  \citet{Bauss2023gitinstances}. }
\end{table}

\paragraph{Results of Generalized Assignment Problems}
Preliminary numerical tests on a different test set show, that the node selection with hypervolume of the search zone box NS(HSZ) works better than the one with the local hypervolume gap NS(LHG). This node selection strategy has a significant impact on the number of explored nodes and on the computational time as well. The number of nodes is reduced by up to $53.4\%$ (Table~\ref{GAP3-18}), while the average run time is reduced by up to $38.9\%$ (Table~\ref{GAP3-18}). 

We can observe, that in the tri-objective instances WST outperforms EC, but in the benchmark instances with four objectives EC outperforms WST. Therefore, the best choice regarding the computation time for instances with three objectives is WST, which reduces the computation time by up to $40.5\%$ (Table~\ref{GAP3-18}). For the instances with four objectives EC seems to be the best choice regarding the computation time with a reduction of up to $38.6\%$ (Table~\ref{GAP4-12}). Regarding average number of explored nodes SLB+TE performs best, since it explores the fewest nodes in the majority of the considered instances. In the best case, the number of nodes is reduced by $64.7\%$ (Table~\ref{GAP4-15}), considering the instance sizes where all problems are solved. If we look at the results in Table~\ref{GAP4-15}, we even see a reduction of nodes by $66.3\%$ and a simultaneous improvement in the number of solved instances. Nevertheless, there are other approaches for this instance size, that solve even more problems.

\paragraph{Summary and General Observations}
In all of the four tested problem classes (KP, UFLP, CFLP and GAP) a significant reduction of the average number of explored nodes and the average total computation time can be realized in nearly every tested adapted branch and bound approach. With a rising number of variables the impact on the performance increases, whereas a rising number of objective functions results in a decrease of the performance. This is due to the general shortcomings of multi-objective branch and bound, struggling with weaker bounds in higher dimensions. In Figure~\ref{fig:performance}, performance profiles of the corresponding best approaches are illustrated. Thereby, the $x$-axis represents the time in seconds and the $y$-axis corresponds to the proportion of solved instances.

In three of the four considered problem classes, preliminary tests have shown, that the dynamic node selection strategy based on the search zone box works better in terms of computation time compared to the local hypervolume gap strategy. Only for knapsack problems the local hypervolume gap is used. The stronger combinatorial structure of UFLP, CFLP and GAP (compared to KP) leads to more complex lower bound sets, which are computationally difficult to handle in the local hypervolume gap strategy. 

Unfortunately, there is no clear winner that performs best on all considered benchmark instances, but we can observe some tendencies. Regarding the number of nodes, SLB+TE seems to be a good choice. For every considered problem class, SLB+TE performs comparatively good and is even the best choice in some of them. This approach can also improve the number of solved instances in a few cases. But, due to the high amount of solved integer programming scalarizations this approach often increases the total computation time, especially for smaller instance sizes. Regarding the total computation time, WST and EC seem to be good choices. Both approaches are the best performing methods in some of the problem classes and perform also comparatively good in the other ones. 

Nevertheless, by just using the corresponding dynamic node selection strategy, remarkable improvements are achieved w.r.t.\ both --- the number of considered branch and bound nodes and the computation time. 

\begin{figure}[htbp!]
	\subfloat[Performance of tri-objective knapsack problems \label{Fig:KP3res}]{
		\begin{tikzpicture}[scale =0.499]
			\begin{axis}[
				width=\linewidth, 
				grid=major, 
				grid style={dashed,gray!30},
				y label style={at={(axis description cs:-0.08,.5)},anchor=south},
				xlabel= time (s), 
				ylabel=proportion of instances solved,
				legend style={at={(0,1)},anchor=north west},
				x tick label style={rotate=90,anchor=east}
				]
				\addplot[color = RoyalBlue, const plot]
				table[x=BB,y=proportion,col sep=comma] {knapsack3.csv}; 
				\addplot[color = Maroon,const plot]
				table[x=EC+TE,y=proportion,col sep=comma] {knapsack3.csv}; 
				\addplot[color = PineGreen, const plot]
				table[x=SLB+TE,y=proportion,col sep=comma] {knapsack3.csv}; 
				\legend{BB,EC+TE,SLB+TE}
				
			\end{axis}
		\end{tikzpicture}
	}
\subfloat[Performance of 4-objective knapsack problems \label{Fig:KP4sres}]{
	\begin{tikzpicture}[scale =0.499]
		\begin{axis}[
			width=\linewidth, 
			grid=major, 
			grid style={dashed,gray!30},
			y label style={at={(axis description cs:-0.08,.5)},anchor=south},
			xlabel= time (s), 
			ylabel=proportion of instances solved,
			legend style={at={(0,1)},anchor=north west},
			x tick label style={rotate=90,anchor=east}
			]
			\addplot[color = RoyalBlue, const plot]
			table[x=BB,y=proportion,col sep=comma] {knapsack4.csv}; 
			\addplot[color = Maroon,const plot]
			table[x=EC+TE,y=proportion,col sep=comma] {knapsack4.csv}; 
			\addplot[color = PineGreen, const plot]
			table[x=SLB+TE,y=proportion,col sep=comma] {knapsack4.csv}; 
			\legend{BB,EC+TE,SLB+TE}
			
		\end{axis}
	\end{tikzpicture}
}

	\subfloat[Performance of tri-objective uncapacitated facility location problems \label{Fig:Uflp3res}]{
		\begin{tikzpicture}[scale =0.499]
			\begin{axis}[
				width=\linewidth, 
				grid=major, 
				grid style={dashed,gray!30},
				y label style={at={(axis description cs:-0.08,.5)},anchor=south},
				xlabel= time (s), 
				ylabel=proportion of instances solved,
				legend style={at={(0,1)},anchor=north west},
				x tick label style={rotate=90,anchor=east}
				]
				\addplot[color = RoyalBlue, const plot]
				table[x=BB,y=proportion,col sep=comma] {uflp3.csv}; 
				\addplot[color = Maroon,const plot]
				table[x=WST,y=proportion,col sep=comma] {uflp3.csv}; 
				\addplot[color = PineGreen, const plot]
				table[x=WST+TE,y=proportion,col sep=comma] {uflp3.csv}; 
				\legend{BB,WST,WST+TE}
				
			\end{axis}
		\end{tikzpicture}
	}
\subfloat[Performance of 4-objective uncapacitated facility location problems \label{Fig:Uflp4res}]{
	\begin{tikzpicture}[scale =0.499]
		\begin{axis}[
			width=\linewidth, 
			grid=major, 
			grid style={dashed,gray!30},
			y label style={at={(axis description cs:-0.08,.5)},anchor=south},
			xlabel= time (s), 
			ylabel=proportion of instances solved,
			legend style={at={(0,1)},anchor=north west},
			x tick label style={rotate=90,anchor=east}
			]
			\addplot[color = RoyalBlue, const plot]
			table[x=BB,y=proportion,col sep=comma] {uflp4.csv}; 
			\addplot[color = Maroon,const plot]
			table[x=WST,y=proportion,col sep=comma] {uflp4.csv}; 
			\addplot[color = PineGreen, const plot]
			table[x=WST+TE,y=proportion,col sep=comma] {uflp4.csv}; 
			\legend{BB,WST,WST+TE}
			
		\end{axis}
	\end{tikzpicture}
}

\subfloat[Performance of tri-objective generalized assignment problems \label{Fig:Gap3res}]{
	\begin{tikzpicture}[scale =0.499]
		\begin{axis}[
			width=\linewidth, 
			grid=major, 
			grid style={dashed,gray!30},
			y label style={at={(axis description cs:-0.08,.5)},anchor=south},
			xlabel= time (s), 
			ylabel=proportion of instances solved,
			legend style={at={(0,1)},anchor=north west},
			x tick label style={rotate=90,anchor=east}
			]
			\addplot[color = RoyalBlue, const plot]
			table[x=BB,y=proportion,col sep=comma] {gap3.csv}; 
			\addplot[color = Maroon,const plot]
			table[x=WST,y=proportion,col sep=comma] {gap3.csv}; 
			\addplot[color = PineGreen, const plot]
			table[x=SLB+TE,y=proportion,col sep=comma] {gap3.csv}; 
			\legend{BB,WST,SLB+TE}
			
		\end{axis}
	\end{tikzpicture}
}
\subfloat[Performance of 4-objective generalized assignment problems \label{Fig:Gap4res}]{
	\begin{tikzpicture}[scale =0.499]
		\begin{axis}[
			width=\linewidth, 
			grid=major, 
			grid style={dashed,gray!30},
			y label style={at={(axis description cs:-0.08,.5)},anchor=south},
			xlabel= time (s), 
			ylabel=proportion of instances solved,
			legend style={at={(0,1)},anchor=north west},
			x tick label style={rotate=90,anchor=east}
			]
			\addplot[color = RoyalBlue, const plot]
			table[x=BB,y=proportion,col sep=comma] {gap4.csv}; 
			\addplot[color = Maroon,const plot]
			table[x=EC,y=proportion,col sep=comma] {gap4.csv}; 
			\addplot[color = PineGreen, const plot]
			table[x=SLB+TE,y=proportion,col sep=comma] {gap4.csv}; 
			\legend{BB,EC,SLB+TE}
			
		\end{axis}
	\end{tikzpicture}
}

\caption{Performance profiles of different problem classes.}
\label{fig:performance}
\end{figure}

\section{Conclusion and Outlook}\label{sec:dis}
In this paper we propose different adaptive improvements for multi-objective branch and bound frameworks using objective-space information to partially overcome its structural difficulties. We propose new dynamic node selection strategies, based on the gap between lower and upper bound set, which improve the number of explored nodes and the total computation time. Furthermore, we use  objective space information, gained by solving scalarized (sub)problems to integer optimality, to improve the lower and upper bound set. Additionally, simple lower bounds set are used adaptively to omit the computation of the complete lower bound set. The numerical results show the positive impact on different problem classes regarding the number of explored nodes and the computational time. 

This paper shows, that the order in which the nodes are explored, has a significant impact on the performance of multi-objective branch and bound algorithms. Since the order depends on the node selection strategy and the branching rule, it might be promising to include different combinations of those components. Especially the combination of different dynamic, problem dependent strategies could result in further improvements.

\paragraph{Acknowledgments}
Partially funded by the by Deutsche Forschungsgemeinschaft (DFG, German Research Foundation), project number 441310140.


\end{document}